\documentclass[aos]{imsart}
\usepackage{graphicx,amsthm,amsmath,cglw,pstricks}
\renewcommand{\E}{{\sf E}}                          
\renewcommand{\R}{{\mathbb{R}}}

\begin{document}

\begin{frontmatter}

% "Title of the paper"
\title{Adaptive Confidence Bands}
\runtitle{Confidence Bands}

% indicate corresponding author with \corref{}
% \author{\fnms{John} \snm{Smith}\corref{}\ead[label=e1]{smith@foo.com}\thanksref{t1}}
% \thankstext{t1}{Thanks to somebody} 
% \address{line 1\\ line 2\\ printead{e1}}
% \affiliation{Some University}

\author{\fnms{Christopher} \snm{Genovese}\ead[label=e1]{genovese@stat.cmu.edu}}
\address{Department of Statistics\\ Carnegie Mellon University\\ \printead{e1}}
\and
\author{\fnms{Larry} \snm{Wasserman}\ead[label=e2]{larry@stat.cmu.edu}}
\address{Department of Statistics\\ Carnegie Mellon University\\ \printead{e2}}
\affiliation{Carnegie Mellon University}

\runauthor{Genovese and Wasserman}

\begin{abstract}
We show that there do not exist adaptive confidence bands for curve
estimation except under very restrictive assumptions.  We propose
instead to construct adaptive bands that cover a surrogate function
$f^\star$ which is close to, but simpler than, $f$.  The surrogate
captures the significant features in $f$. We establish lower bounds
on the width for any confidence band for $f^\star$ and construct
a procedure that comes within a small constant factor
of attaining the lower bound for finite-samples.

\end{abstract}


%\begin{keyword}
%\kwd{}
%\kwd{}
%\end{keyword}

\end{frontmatter}

%\baselineskip=12pt
%
%\begin{center}
%{\large {\bf Adaptive Confidence Bands}}
%\end{center}
%
%\begin{center}
%{\sc By} {\sc Christopher R. Genovese}\footnote{
%Research supported by NSF Grant SES 9866147.}
%{\sc and} {\sc Larry Wasserman}\footnote{
%Research supported by NIH Grant R01-CA54852-07,
%NIH grant number MH57881,
%NSF Grant DMS-98-03433 and
%NSF Grant DMS-0104016.}\\
%\medskip
%\emph{Department of Statistics}\\
%\emph{Carnegie Mellon University}\\
%\end{center}
%
%
%\begin{center}
%{\large{\gray{\sf \today}}}
%\end{center}
%
%
%\begin{quote}
%  
%
%{\sc Key words and phrases:}
%Confidence sets, confidence bands,
%nonparametric regression.
%\end{quote}
%
%
%\baselineskip=16pt
%

\section{Introduction}

\subsection{Motivation}

Let $(x_1,Y_1), \ldots, (x_n,Y_n)$ 
be observations from the nonparametric regression model
\begin{equation}
Y_i = f(x_i) + \sigma\,\epsilon_i
\end{equation}
where $\epsilon_i \sim N(0,1)$,
$x_i \in (0,1)$,
and $f$ is assumed to lie in some infinite-dimensional class of functions $\cH$.
We are interested in constructing confidence bands $(L, U)$ for $f$.
Ideally these bands should satisfy
\begin{equation}\label{eq::true-cov}
\P_f{L \le f \le U} = 1-\alpha \ \ \ \ \ \
\mbox{for all}\ f\in \cH
\end{equation}
where
${L \le f \le U}$ means that
$L(x) \le f(x) \le U(x)$
for all $x\in\cX$,
where $\cX$ is some subset of $(0,1)$
such as
$\cX=\{x\}, \cX = \{x_1, \ldots, x_n\}$ or
$\cX = (0,1)$.
Throughout this paper, we take
$\cX = \{x_1, \ldots, x_n\}$ but this particular choice
is not crucial in what follows.

Attaining (\ref{eq::true-cov}) is difficult and hence
it is common to settle for
pointwise asymptotic coverage:
\begin{equation}\label{eq::pointwise}
\liminf_{n\to\infty}\P_f{L \le f \le U} \geq 1-\alpha \ \ \ \ \ \
\mbox{for all}\ f\in \cH.
\end{equation}
``Pointwise'' refers to the fact that the asymptotic limit
is taken for each fixed $f$ rather than uniformly over $f\in\cH$.
Papers on pointwise asymptotic methods include 
Claeskens and Van Keilegom (2003),
Eubank and Speckman (1993),
H\"{a}rdle and Marron (1991),
Hall and Titterington (1988),
H\"{a}rdle and Bowman (1988),
Neumann and Polzehl (1998), 
and Xia (1998).

Achieving even pointwise asymptotic coverage
is nontrivial due to the presence of bias.
If $\hat{f}(x)$ is an estimator with mean
$\overline{f}(x)$ and standard deviation
$s(x)$ then
$$
\frac{\hat{f}(x) - f(x)}{s(x)} = 
\frac{\hat{f}(x) - \overline{f}(x)}{s(x)} +
\frac{{\rm bias}(x)}{\sqrt{\rm variance(x)}}.
$$
The first term typically satisifes a central
limit theorem but the second term 
does not vanish
even asymptotically if the bias and variance are
balanced.
For discussions on this point, see the papers
referenced above as well as
Ruppert, Wand, and Carroll (2003) and Sun and Loader (1994).

Pointwise asymptotic bands are not uniform, that is, they do not
control
\begin{equation}
\inf_{f\in \cH}\P_f{L \le f \le U}.
\end{equation}
The sample size
$n(f)$ required 
for the true coverage to approximate the nominal coverage,
depends on the unknown function $f$.

The aim of this paper is to attain uniform coverage over $\cH$.
We say that
$B=(L,U)$ has \emph{uniform coverage} if
\begin{equation}\label{eq::first}
\inf_{f\in\cH}\P_f{L \le f \le U}  \ge 1-\alpha .
\end{equation}
Starting in Section \ref{sec::projections},
we will insist on coverage over
$\cH = \{{\rm all\  functions}\}$.

The bound  in (\ref{eq::first}) can be achieved trivially
using Bonferroni bands.
Set $\ell_i= Y_i - c_n \sigma$ and
$u_i= Y_i + c_n \sigma$,
where $c_n = \Phi^{-1}(1-\alpha/2n)$ and
$\Phi$ is the standard Normal {\sc cdf}.
Yet this band is unsatisfactory for several reasons: 
\begin{enumerate}
\item The width of the band grows with sample size.
\item The band is centered on a poor estimator of the unknown function.
\item The width of the band is independent of the data and hence cannot adapt
to the smoothness of the unknown function.
\end{enumerate}
Problems (1) and (2) are easily remedied by using standard smoothing methods.
But the results of Low (1997) suggest that (3) is an inevitable consequence of uniform coverage.

The smoother the functions in $\cH$, the smaller the width necessary to achieve uniform coverage.
Suppose that $\cF \subset \cH$ contains the ``smooth'' functions in $\cH$ and that $\cH - \cF$ is nonempty.
Uniform coverage over $\cH$ requires that
the width of fixed-width bands be driven by the ``rough'' functions in $\cH - \cF$; the width will thus be large even if $f\in\cF$.
Ideally, our procedure would adjust automatically to produce narrower bands when the function 
is smooth ($f\in\cF$) and wider bands when the function is rough ($f\not\in\cF$),
but to do that, the width must be determined from the data.
Low showed that for density estimation at a single point,
fixed-width confidence intervals perform as well as random length intervals;
that is, the data do not help reduce the width of the bands for smoother functions.
In Section \ref{sec::failure}, 
we extend Low's result to nonparametric regression
and show that the phenomenon is quite general.
Without restrictive assumptions, confidence bands cannot adapt.

These results mean that the width of uniform confidence bands is determined by the
greatest roughness we are willing to assume.
Because the typical assumptions about $\cH$ in the nonparametric regression problem are
loosely held and difficult to check, 
the result is that the confidence band widths are essentially arbitrary.  
This is not satisfactory in practice.

The contrast with $L^2$ confidence balls is noteworthy.
$L^2$ confidence sets
have been studied by
Li (1999), Juditsky and Lambert-Lacroix (2002), 
Beran and D\"{u}mbgen (1998), Genovese and Wasserman (2004),
Baraud (2004), Hoffman and Lepski (2003),
Cai and Low (2004), and Robins and van der Vaart (2004).
Let
\begin{equation}
B = \Biggl\{{ f}\in \mathbb{R}^n:\ \frac{1}{n}\sum_{i=1}^n({ f}_i-\hat{{ f}}_i)^2 \le R_n^2\Biggr\}
\end{equation}
for some $\hat{f}$
and suppose that
\begin{equation}
\inf_{{ f}\in\mathbb{R}^n}\P_f{{ f}\in B}\ge 1-\alpha .
\end{equation}
Then
\begin{equation}\label{eq::smaller}
\inf_{{ f}\in\mathbb{R}^n}\E_f(R_n) \ge \frac{C_1}{n^{1/4}},
\ \ \ \ \mbox{and}\ \ \ \ 
\sup_{{ f}\in\mathbb{R}^n}\E_f(R_n) \ge C_2
\end{equation}
where $C_1$ and $C_2$ are positive constants.
Moreover, there exist confidence sets that
achieve the faster $n^{-1/4}$ rate at some points
in $\mathbb{R}^n$.
Because fixed-radius confidence sets necessarily have radius of size O(1),
the supremum in (\ref{eq::smaller}) implies such confidence sets
must have random radii.
We can construct random-radius confidence balls that improve on fixed-radius confidence sets,
for example, by obtaining a smaller radius for subsets of smoother functions $f$.
$L^2$ confidence balls can therefore adapt to the unknown smoothness of $f$.
Unfortunately, confidence balls can be difficult to work with in high dimensions (large $n$)
and tend to constrain many features of interest rather poorly,
for which reasons confidence bands are often desired.

It is also interesting to compare the adaptivity results for estimation and inference.
Estimators exist (e.g., Donoho et al. 1995) that can adapt to unknown smoothness, achieving
near optimal rates of convergence over a broad scale of spaces.
But since confidence bands cannot adapt, the minimum width bands that achieve uniform coverage
over the same scale of spaces have width $O(1)$, overwhelming the differences among reasonable estimators.
We are left knowing that we are close to the true function but being unable to demonstrate it
inferentially.

The message we take from the nonadaptivity results
in Low (1987) and Section \ref{sec::failure} of this paper
is that the problem of constructing confidence
bands for $f$ over nonparametric classes is simply too difficult
under the usual definition of coverage.
Instead, we introduce a slightly weaker notion -- surrogate coverage -- 
under which it is possible to obtain adaptive bands while allowing sharp
inferences about the main features of $f$. %% NOTE: don't like ``main'' here...

\subsection{Surrogates}

\begin{figure}
\hspace{1cm}
\includegraphics[width=5in]{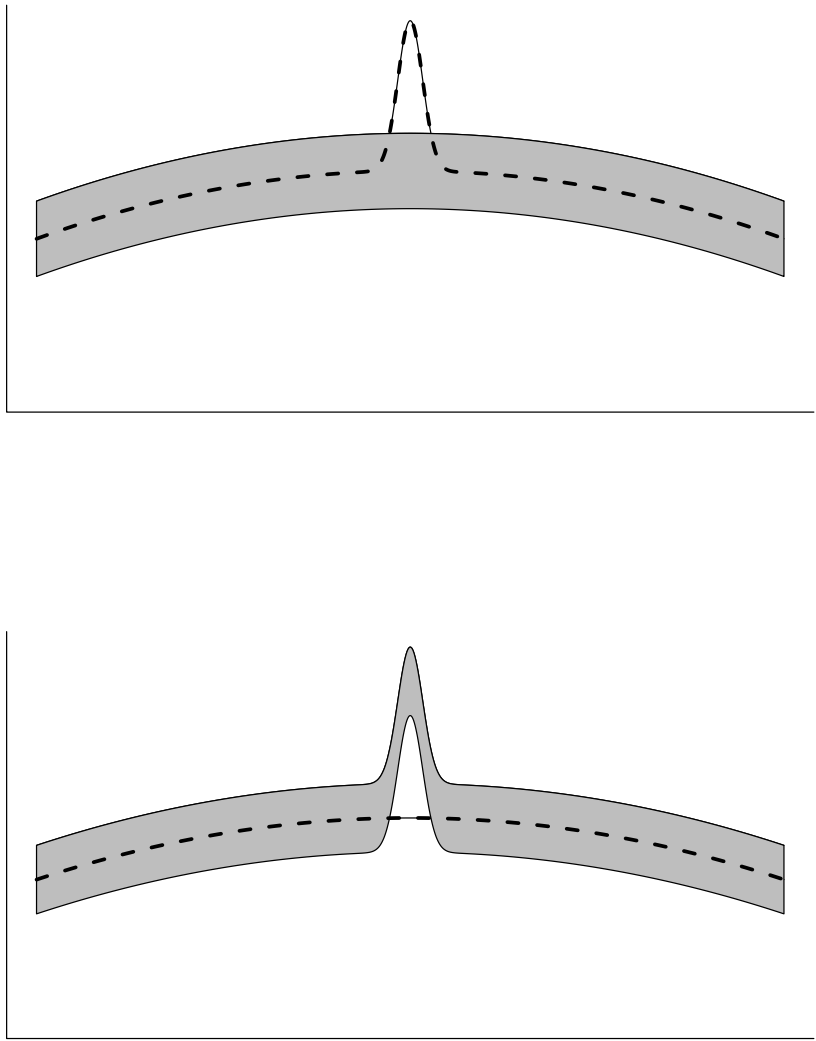}
\caption{The top plot shows a conservative failure: the only place where $f$ is not
contained in the band is when the bands are smoother than the truth.
The bottom plot
shows a liberal failure: the only place where $f$ is not
contained in the band is when the bands are less smooth than the truth.
The usual notion of coverage treats these failures equally.}
\label{fig::bandfig}
\end{figure}

Figure \ref{fig::bandfig} shows two situations where a band fails to
capture the true function.
The top plot shows a conservative failure: the only place where $f$ is not
contained in the band is when the bands are smoother than the truth.
The bottom plot
shows a liberal failure: the only place where $f$ is not
contained in the band is when the bands are less smooth than the truth.
The usual notion of coverage treats these failures equally.
Yet, in some sense, the second error is more serious than the first
since the bands overstate the complexity.

\begin{figure}
\hspace{1cm}
\includegraphics[width=5in]{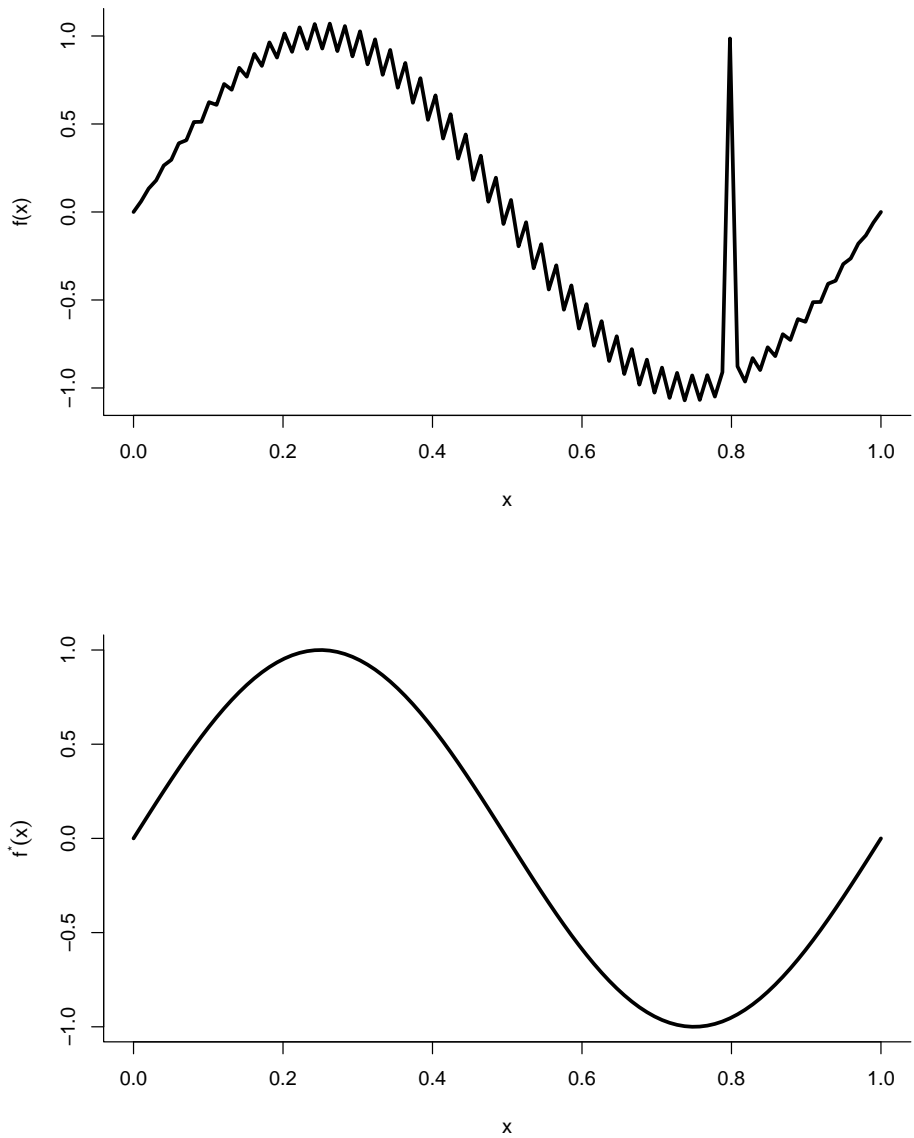}
\caption{The top plot shows a complicated function $f$.
The bottom shows a surrogate $f^\star$ which is simpler than $f$
but retains the main, estimable features of $f$.
Adaptation is possible if we cover $f^\star$ instead of $f$.}
\label{fig::surrogate}
\end{figure}

We are thus led to a different approach
that treats conservative errors and liberal errors differently.
The basic idea 
is to find a function $f^\star$ that is simpler than $f$
as in Figure \ref{fig::surrogate}.
We then require that
\begin{equation}
\P_f{ L \le f \le U\ {\rm or}\ L \le f^\star \le U} \ge 1-\alpha,\ \ \ 
{\rm for \ all\ functions\ }f.
\end{equation}
More generally, we will define a finite set of
surrogates $F^\star \equiv F^*(f) = \{f, f_1^*,\ldots, f_m^*\}$ 
and require that a surrogate confidence band $(L,U)$ satisfy 
\begin{equation}
\inf_f \P_f{ L \le g \le U\ \ {\rm for\ some\ }g\in F^\star} \ge 1-\alpha.
\end{equation}
We will also consider bands that are adaptive in the following sense:
if $f$ lies in some subspace $\cF$,
then with high probability $\norm{U-L}_\infty \le w(\cF)$,
where $w(\cF)$ is the best width of a uniformly valid confidence band (under the usual definition of coverage)
based on the a priori knowledge that $f\in\cF$.
Among possible surrogates,
a surrogate will be optimal if it admits a valid, adaptive procedure
and the set $\Set{f\in\cF\st F^*(f) = \{f\}}$ is as large as possible.

\subsection{Summary of Results}

In Section \ref{sec::failure},
we show that Low's result on density estimation
holds in regression as well.
Fixed width bands do as well as random width bands,
thus ruling out adaptivity.
We show this when $\cH$ is the set of all functions
and when $\cH$ is a ball in a Lipschitz, 
Sobolev, or Besov space.

Section \ref{sec::projections} gives our main results.
Theorem \ref{thm::main} 
establishes lower bounds
on the width for any valid surrogte confidence band.
Let $\cF$ be a subspace of dimension $d$ in $\R^n$.
The functions that prevent adaptation are those that are
close to $\cF$ in $L^2$ but far in $L^\infty$.
Loosely speaking, such functions are close to $\cF$ except for isolated, spiky features.
If $||f- \Pi f||_2 < \epsilon_2$ and $||f- \Pi f||_\infty > \epsilon_\infty$, 
for tuning constants $\epsilon_2, \epsilon_\infty$,
define the surrogate $f^\star$ to be the projection of $f$ onto $\cF$, $\Pi f$.
Otherwise, define $f^\star = f$.
We show that if 
$\P_f{\norm{U - L}_\infty < w} \ge 1-\gamma$ for all $f\in \cF$,
then
\begin{equation}\label{eq::summaryw}
w \ge \max\left( w_\cF(\alpha,\gamma,\sigma), v(\epsilon_2,\epsilon_\infty, n, d, \alpha, \gamma, \sigma) \right),
\end{equation}
where $w_\cF$ is the minimum width for a uniform confidence band knowing a priori that $f\in\cF$
and
$v(\epsilon_2,\epsilon_\infty, n, d, \alpha, \gamma)$
is described later.

Corollary \ref{cor::optimality} shows that for proper choice of $\epsilon_2$ and $\epsilon_\infty$,
the $v$ term in the previous equation can be made smaller than $w_\cF$.
Figure \ref{fig::surrogate-pic} represents the functions involved;
the gray shaded area are those functions that are replaced by surrogates in the coverage statement,
denoted later by $\cS(\epsilon_2,\epsilon_\infty)$.
These are the functions that are both hard to distinguish from $\cF$ (because they are close to it)
and hard to cover (because they are ``spiky'').
The optimal choice of $\epsilon_2$ and $\epsilon_\infty$ minimizes the volume of this set
while making the right hand side in inequality (\ref{eq::summaryw}) equal to $w_\cF$.
Put another way,
the richest model that permits adaptive confidence bands under the usual notion of coverage
is $\cF =   \mathbb{R}^n - \cS(\epsilon_2,\epsilon_\infty)$.

Theorem \ref{thm::achieve-multi}
gives a procedure that comes within a factor of 2
of attaining the lower bound for finite-samples.
The procedure conducts goodness of fit tests for subspaces
and constructs bands centered on the estimator of the
lowest dimensional nonrejected subspace.
Such a procedure actually reflects common practice.
It is not uncommon to fit a model, check the fit, and
if the model does not fit then we fit a more complex model.
In this sense, we view our results as providing a rigorous basis
for common practice.
It is known that pretesting followed by inference
does not lead to valid inferences for $f$ (Leeb and and P\"otscher, 2005). 
But if we cant accept that sometimes we cover a surrogate $f^\star$ rather than $f$,
then validity is restored.

These results are proved in Section \ref{sec::proofs}.

\subsection{Related Work}

The idea of estimating the detectable part of $f$ is present, at least
implicitly, in other approaches.  Davies and Kovac (2001) separate the
data into a simple piece plus a noise piece which is
similar in spirit to our approach.  
Another related idea is
scale-space inference due to Chaudhuri and Marron (2000) who focus
on inference for all smoothed versions of $f$ rather than $f$ itself.
Also related is the idea of oversmoothing
as described in Terrell (1990) and
Terrell and Scott (1985).
Terrell argues that 
``By using the most smoothing that is compatible 
with the scale of the problem, we tend to eliminate
accidental features.''
The idea of one-sided inference
in Donoho (1988) has a similar spirit.
Here, one constructs confidence intervals
of the form $[L,\infty)$
for functionals such as the number of modes
of a density.
Bickel and Ritov (2000)
make what they call a ``radical proposal''
to `` ... determine how much bias can be tolerated without
[interesting] features being obscured.''
We view our approach as a way of implementing their suggestion.
Another related idea is contained in Donoho (1995) who showed that
if $\hat{f}$ is the soft threshold estimator of a function and
$f(x) = \sum_j \theta_j \psi_j(x)$ is an expansion in an unconditional basis,
then
$\P_f{\hat{f}\preceq f} \ge 1-\alpha$
where
$\hat{f}=\sum_j \hat\theta_j \psi_j$ and
$\hat{f}\preceq f$ means that
$|\hat\theta_j| \le |\theta_j|$ for all $j$.
Finally, we remind the reader that there is a plethora of work on
adaptative estimation; see, for example, Cai and Low (2004)
and references therein.

\subsection{Notation}

If $L$ and $U$ are random functions on $\cX = \Set{x_1,\ldots,x_n}$ such that $L \le U$,
we define $B = (L,U)$ to be the (random) set of all functions $g$ on $\cX$ for which $L \le g \le U$.
We call $B$ (or equivalently, the pair $L,U$) a band;
the band covers a function $f$ if $f\in B$ (or equivalently, if $L \le f \le U$).
Define its width to be the random variable
\begin{equation}
W = \norm{U - L}_\infty = \max_{1\le i\le n} (U(x_i) - L(x_i)). 
\end{equation}

Because we are constructing bands on $\cX = \{x_1,\ldots,x_n\}$, we most often
refer to functions in terms of their evaluations $f = (f(x_1),\ldots,f(x_n)) \in \R^n$.
When we need to refer to a space of functions to which $f$ belongs, we use a $\;\tilde{}\;$ to
denote the function space and no $\;\tilde{}\;$ to denote the vector space of evaluations. 
Thus, if $\tilde\cA$ is the space of all functions, then $\cA = \R^n$.
In both cases, we use the same symbol for the function and let the meaning be clear from context;
for example, $f\in\tilde\cA$ is the function and $f\in\cA$ is the vector $(f(x_1),\ldots,f(x_n))$.
Define the following norms on $\mathbb{R}^n$:
\begin{eqnarray*}
||f|| &=& ||f||_2 = \sqrt{\frac{1}{n}\sum_{i=1}^n f_i^2}\\
||f||_\infty &=& \max_i |f_i|.
\end{eqnarray*}
We use $\langle \cdot,\cdot\rangle$ to denote the inner product $\langle f,g\rangle = \frac{1}{n}\sum_{i=1}^n f_i g_i$
corresponding to $\norm{\cdot}$.

If $\cF$ is a subspace of $\R^n$, we define $\Pi_\cF$ to be the Euclidean projection onto $\cF$,
using just $\Pi$ if the subspace is clear from context.
We use
\begin{equation}\label{eq::eibasis}
e_i = (\underbrace{0,\ldots,0}_{i-1\ {\rm times}},1, 
\underbrace{0,\ldots, 0}_{n-i\ {\rm times}})^T
\end{equation}
to denote the standard basis on $\R^n$.

If $F_\theta$ is a family of {\sc cdf}s indexed by $\theta$,
we write $F_\theta^{-1}(\alpha)$ to denote the lower-tail $\alpha$-quantile of $F_\theta$.
For the standard normal distribution, however,
we use $z_\alpha$ to denote the upper-tail $\alpha$-quantile,
and we denote the {\sc cdf} and {\sc pdf}, respectively, by $\Phi$ and $\phi$.

Throughout the paper we assume that $\sigma$ is a known constant;
in some cases we simply set $\sigma =1$.
But see Remark \ref{remark::sigma} about the unknown $\sigma$ case.

%%Note: more common notation that should go here?? 
%%Otherwise, is this really necessary? Seems ok but this could be absorbed. 

\section{Nonadaptivity of Bands}\label{sec::failure}

In this section we construct
lower bounds on the width of valid
confidence bands
analagous to
(\ref{eq::smaller})
and we show that the lower bound is achieved by fixed-width bands.

Low (1997)
considered estimating a density $f$
in the class
$$
\cF(a,k,M)=\Biggl\{ f:\ f\ge 0,\,\int f=1,\, f(x_0)\le a,\, ||f^{(k)}(x)||_\infty \le M \Biggr\}.
$$
He shows that if $C_n$ is a confidence interval for $f(0)$, that is,
$$
\inf_{f\in \cF(a,k,M)} \P_f{f(0)\in C_n} \ge 1-\alpha ,
$$
then, for every $\epsilon >0$,
there exists $N = N(\epsilon, M)$ and $c>0$ such that, for all $n \ge N$,
\begin{equation}
\E_f ({\rm length}(C_n)) \ge c \, n^{-k/(2k+1)}
\end{equation}
for all $f\in \cF(a,k,M)$ such that $f(0) > \epsilon$.
Moreover, there exists a fixed-width confidence interval $C_n$
and a constant $c_1$ such that
$\E_f ({\rm length}(C_n)) \le c_1 n^{-k/(2k+1)}$
for all $f\in \cF(a,k,M)$.
Thus, the data play no role in constructing a rate-optimal band,
except in determining the center of the interval.

For example, if we use kernel density estimation,
we could construct an optimal bandwidth $h = h(n,k)$
depending only on $n$ and $k$ -- but not the data --
and construct the interval from that kernel
estimator. This makes the interval highly dependent
on the minimal amount of smoothness $k$ that is assumed.
And it rules out the usual data-dependent bandwidth methods
such as cross-validation.

Now return to the regression model
\begin{equation}\label{eq::normal-means-model}
Y_i = f_i + \sigma \, \epsilon_i,\ \ \ i=1, \ldots, n,
\end{equation}
where $\epsilon_1$, $\ldots$, $\epsilon_n$ are independent,
${\rm Normal} (0,1)$ random variables, and
$f = (f_1, \ldots, f_n) \in \mathbb{R}^n$.

\begin{theorem}
\label{thm::basic}
Let
$B = (L,U)$
be a $1- \alpha$ confidence band
over $\Theta$, where
$0 < \alpha < 1/2$
and let $g\in \Theta$.
Suppose that $\Theta$ contains
a finite set of vectors $\Omega$,
such that:
\begin{enumerate}
\item for every distinct pair
$f,\nu\in\Omega$, we have
$\langle f-g, \nu-g \rangle =0$ and
\item for some $0 < \epsilon < (1/2) - \alpha$,
\begin{equation}
\max_{f\in\Omega}\frac{e^{n ||f-g||^2/\sigma^2}}{|\Omega|} \le  \epsilon^2.
\end{equation}
\end{enumerate}
Then, 
\begin{equation}
\E_g(W)\ge  (1-2\alpha - 2\epsilon) \min_{f\in\Omega}||g-f||_\infty .
\end{equation}
\end{theorem}

We begin with the case where $\Theta = \mathbb{R}^n$.  We will obtain
a lower bound on the width of any confidence band and then show that a
fixed-width procedure attains that width.  The results hinge on
finding a least favorable configuration of mean vectors that are as
far away from each as possible in $L^\infty$ while staying a fixed
distance $\epsilon$ in total-variation distance.

\begin{theorem}
\label{thm::finite-sample-lower-bound}
Let $\cH = \mathbb{R}^n$
and fix $0 < \alpha < 1/2$.
Let
$B = (L,U)$
be a $1- \alpha$ confidence band over $\cH$.
Then, for every $0 < \epsilon < (1/2) - \alpha$,
\begin{equation}\label{eq::rnbound}
\inf_{f\in\mathbb{R}^n}\E_f(W)\ge  (1-2\alpha - 2\epsilon)\sigma \sqrt{\log( n \epsilon^2 )}.
\end{equation}
The bound is achieved (up to constants) by the fixed-width Bonferroni bands:
$$
\ell_i = Y_i - \sigma z_{\alpha/n},\ u_i = Y_i + \sigma z_{\alpha/n}.
$$
\end{theorem}

\begin{theorem}[Lipshschitz Balls]
\label{thm::finite-sample-lip}
Define $x_i = i/n$ for $1 \le i \le n$. Let 
\begin{eqnarray}
\tilde\cH(L) &=& \Biggl\{ f:\ |f(x) - f(y)| \le L |x-y| ,\ \ \, x,y \in [0,1]\Biggr\},\\
\noalign{\noindent be a ball in Lipschitz space, and let}
\cH(L) &=& \{ (f(x_1), \ldots, f(x_n)):\ f\in \tilde\cH(L)\}
\end{eqnarray}
be the vector of evaluations on $\cX$
Fix $0 < \alpha < 1/2$ and let $B = (L,U)$
be a $1- \alpha$ confidence band
over $\cH(L)$.
Then, for every $0 < \epsilon < (1/2) - \alpha$,
\begin{equation}
\inf_{f\in\cH(L)}\E_f(W)  \ge  a_n
\end{equation}
where
\begin{eqnarray*}
a_n &=& \left(\frac{\log n}{n}\right)^{1/3}\times
\left(\frac{L\sigma^2}{2}\right)^{1/3}\\
&& \hspace{-1cm}\times
\left( 1 +  \frac{3\log (1+\epsilon^2)}{\log n} +
 \frac{2 \log (L/(2\sigma))}{\log n} -
 \frac{\log\left(\frac{1}{3}\log n + \log(1+\epsilon^2)+ \frac{2}{3}\log(L/(2\sigma))\right)}
      {\log n} \right).
\end{eqnarray*}
The lower bound is achieved (up to logarithmic factors) by a fixed-width procedure.
\end{theorem}

\begin{theorem}[Sobolev Balls]
\label{thm::sobolev}
Let $\tilde\cH(p,c)$ be a Sobolev ball of order $p$
and radius $c$
and let $B = (L,U)$
be a $1- \alpha$ confidence band over $\cH(p,c)$.
For every $0 < \epsilon < (1/2) - \alpha$, 
for every $\delta >0$,
and all large $n$,
\begin{equation}
\inf_{F\in\cH(p,c-\delta)}
\E_F(W)\ge  (1-2\alpha - 2\epsilon) \left(\frac{c_n}{n^{p/(2p+1)}}\right)
\end{equation}
for some $c_n$ that increases at most logarithmically.
The bound is achieved (up to logarithmic factors) by a fixed-width band procedure.
\end{theorem}

\begin{theorem}[Besov Balls] \label{thm::besov}
Let $\tilde\cH(p,q,\xi,c)$ be ball of size $c$ in the Besov space $B_{p,q}^\xi$
and
et $B = (L,U)$ be a $1- \alpha$ confidence band over $\cH(p,q,\xi,c)$.
For every $0 < \epsilon < (1/2) - \alpha$, and every $\delta>0$,
\begin{equation}
\inf_{f\in\cH(p,q,\xi,c-\delta)}\E_f(W)\ge c_n (1-2\alpha - 2\epsilon) n^{-1/(1/p - \xi - 1/2)}.
\end{equation}
The bound is achieved (up to logarithmic factors) by a fixed-width procedure.
\end{theorem}

\section{Adaptive Bands} \label{sec::projections}

Let $\{{\cal F}_T: T\in \cT\}$ be a scale of 
linear subspaces.
Let $w_T$ denote the smallest width
of any confidence band when it is known that $f\in\cF_T$
(defined more precisely below).
We would like to define an approporiate surrogate
and a procedure that gets as close as possible to
the target width $w_T$ when $f\in\cF_T$.
To clarify the ideas, subsection \ref{sec::singlesubspace}
develops our results in the special case where the subspaces
are $\Set{\cF, \R^n}$ for a fixed $\cF$ of dimension $d < n$.
Subsection \ref{sec::nestedspaces} handles the more general
case of a sequence of nested subspaces. 

\subsection{Preliminaries}

We begin by defining several quantities that will be used throughout.
Let $\tau(\epsilon)$ denote the total variation distance
between a $N(0,1)$ and a $N(\epsilon,1)$ distribution.
Thus,
\begin{equation}
\tau(\epsilon) = \Phi(\epsilon/2) - \Phi(-\epsilon/2).
\end{equation}
Then,
$\epsilon \phi(\epsilon/2)\le \tau(\epsilon)\le \epsilon \phi(0)$
and
$\tau(\epsilon)\sim \epsilon \phi(0)$ as $\epsilon\to 0$.

\begin{lemma}
If $P = N(f,\sigma^2 I)$ and $Q = N(g,\sigma^2 I)$
are multivariate Normals with $f,g\in\mathbb{R}^n$ then
\begin{equation}
d_{\rm TV}(P,Q) = \tau \left(\frac{\sqrt{n}||f-g||}{\sigma}\right).
\end{equation}
\end{lemma}

We will need several constants.
For $0 < \alpha < 1$ and $0 < \gamma < 1-2\alpha$ define
\begin{equation} \label{eq::const::kappa}
\kappa(\alpha,\gamma) 
  = \left(2\log(1 + 4(1 - \gamma - 2\alpha)^2)\right)^{1/4}. 
%% OLD VERSION  = (1/6)\sqrt{2\log(1 + 4(1 - \gamma - 2\alpha)^2)}. 
\end{equation}

For $0 < \beta < 1-\xi < 1$ and
integer $m\ge 1$ define
$Q=Q(m,\beta,\xi)$ to be the solution of
\begin{equation}\label{eq::const::Q}
\xi = 1- F_{0,m}(F^{-1}_{Q\sqrt{m},m}(\beta)),
\end{equation}
where $F_{a,d}$ denotes the {\sc cdf} of a $\chi^2$ 
random variable with $d$ degrees of
freedom and noncentrality parameter $a$

\begin{lemma}\label{lemma::universal::Q}
There is a universal constant $\Lambda(\beta,\xi)$
such that
$Q(m,\beta,\xi)\le \Lambda(\beta,\xi)$
for all $m\ge 1$.
For example, $\Lambda(.05,.05) \le 6.25$.
Suppose now that $m = m_n$, $\beta = \beta_n$, and
$\xi = \xi_n$ are all functions of $n$.
As long as
$-\log\beta_n \le \log n$
and $-\log\xi_n \le \sqrt{\log n}$,
then $Q(m_n,\beta_n,\xi_n) = O(\sqrt{\log n})$.
\end{lemma}

Next, define
\begin{equation}\label{eq::const::E}
E(m,\alpha,\gamma)  = \max( Q(m,\alpha,\gamma), 2 \kappa(\alpha,\gamma)),
\end{equation}
for $0 < \alpha < 1$ and $0 < \gamma < 1-2\alpha$.

Finally,
if $\cF$ is a subspace of dimension $d$,
define
\begin{equation}
\Omega_\cF =
\max_{1\le i\le n} \frac{\norm{\Pi_\cF e_i}}{\norm{e_i}},
\end{equation}
where $e_i$ is defined in equation (\ref{eq::eibasis}).
Note that $0 \le \Omega_\cF \le 1$.
The value of $\Omega_\cF$ relates to the geometry of $\cF$ as a hyperplane
embedded in $\R^n$,
as seen through the following results.

\begin{lemma}\label{lemma::min2inInfinity}
Let $\cF$ be a subspace of $\mathbb{R}^n$. Then
\begin{eqnarray}
\min\Biggl\{ \norm{v}:\ v\in\cF, \norm{v}_\infty = \epsilon \Biggr\} &=&
\frac{\epsilon}{\sqrt{n}\Omega_{\cF}}\\
\max\Biggl\{ \norm{v}_\infty:\ v\in\cF, \norm{v} = \epsilon \Biggr\} &=&
\epsilon \sqrt{n}\Omega_{\cF}.
\end{eqnarray}
\end{lemma}

\begin{lemma}
Let $\{\phi_1, \ldots, \phi_d\}$ be orthonormal
vectors with respect to $||\cdot ||$ in $\mathbb{R}^n$
and let $\cF$ be the linear span of these vectors.
Then
\begin{equation}
\Omega_{\cF} = \max_{1\le i \le n} \sqrt{ \frac{\sum_{j=1}^d \phi_{ji}^2}{n}}.
\end{equation}
In particular, if $\max_j\max_i \phi_j(i)\le c$ then
\begin{equation}
\Omega_{\cF} \le c\sqrt{\frac{d}{n}}.
\end{equation}
\end{lemma}

\begin{lemma}
Let 
$\{\phi_1, \ldots, \phi_d\}$
be orthonormal functions on $[0,1]$.
Define
$\cH_j$ to be the linear span of
$\{\phi_1, \ldots, \phi_j\}$.
Let $x_i=i/n$, $i=1, \ldots, n$ and
$\cF_j = \{ f = (h(x_1), \ldots, h(x_n)):\ h\in\cH_j\}$.
Then, 
\begin{equation}
\Omega_{\cF} = \sqrt{ \frac{\sum_{j=1}^d \phi_j^2(x_i)}{n}} + O(1/n).
\end{equation}
In particular, if $\max_j\sup_x \phi_j(x)\le c$ then
\begin{equation}
\Omega_{\cF} \le c\sqrt{\frac{d}{n}} + O(1/n).
\end{equation}
\end{lemma}

In addition, we need the following Lemma first proved, in a related form,
in Baraud (2003). 

\begin{lemma}\label{lemma::baraud-doover}

Let $\cF$ be a subspace of dimension $d$.
Let $0 < \delta < 1 - \xi$ and 
\begin{equation}
\epsilon = \frac{(n-d)^{1/4}}{\sqrt{n}}\,\left(2 \log(1 + 4 \delta^2)\right)^{1/4}. 
\end{equation}
Define $A = \Set{f:\ \norm{f - \Pi_\cF f} > \epsilon}$.
Then,
\begin{equation}
\beta \equiv
   \inf_{\phi_\alpha\in\Phi_\xi}\sup_{f\in A}\P_f{\phi_\xi=0} \ge 1 - \xi - \delta
\end{equation}
where
\begin{equation}
\Phi_\xi = \Biggl\{ \phi_\xi:\ \sup_{f\in\cF}\P_f{\phi_\xi = 0} \le \xi \Biggr\}
\end{equation}
is the set of level $\xi$ tests.

\end{lemma}

\subsection{Single Subspace}\label{sec::singlesubspace}

To begin, we start with a single subspace $\cF$
of dimension $d$.

\begin{definition}
For given $\epsilon_2, \epsilon_\infty >0$,
define the {\em surrogate} $f^\star$ of $f$ by
\begin{eqnarray}
f^\star=
\left\{\begin{array}{ll}
\Pi f & {\rm if\ }||f-\Pi f||_2 \le \epsilon_2\ {\rm and}\ ||f-\Pi f||_\infty > \epsilon_\infty\\
 f & {\rm otherwise.}
\end{array}
\right.
\end{eqnarray}
Define the \emph{surrogate set} of $f$,  $F^*(f) = \Set{ f, f^* }$, which will be a singleton when $f^* = f$.
Define the {\em spoiler set}
$\cS(\epsilon_2,\epsilon_\infty) = \{f\in\R^n:\ f^\star \ne f\}$
and the {\em invariant set}
$\cI(\epsilon_2,\epsilon_\infty) = \{f:\ f^\star = f\}$.
\end{definition}

\begin{figure}
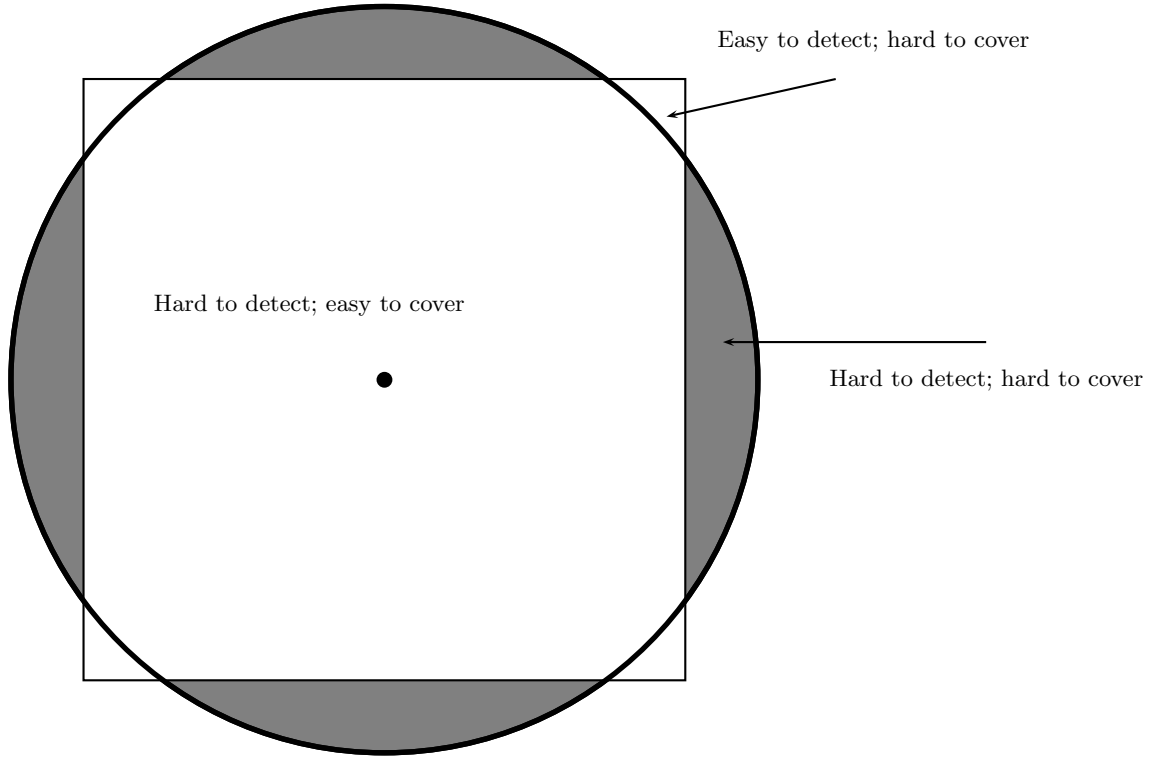

\hglue -2in
\vbox
{
\pspicture(0,3)(16,12)
\pscircle[fillstyle=solid,fillcolor=gray,linewidth=2pt](8,8){5}
\pspolygon[fillstyle=solid,fillcolor=white](12,4)(12,12)(4,12)(4,4)
\pscircle[linewidth=2pt](8,8){5}
\psdots[dotsize=6pt](8,8)
\rput(7,9){Hard to detect; easy to cover}
\rput(16,8){Hard to detect; hard to cover}
\rput(14.5,12.5){Easy to detect; hard to cover}
\psline{->}(16,8.5)(12.5,8.5)
\psline{->}(14,12)(11.75,11.5)
\endpspicture
}
\caption{The dot at the center represents the subspace $\cF$.
The shaded area is the set of spoilers $\cS(\epsilon_2,\epsilon_\infty)$
of vectors for which $f^\star \neq f$. If these vectors were not surrogated,
adaptation is not possible.
The non-shaded area is the invariant set $\cI(\epsilon_2,\epsilon_\infty) =
\{f:\ f^\star = f\}$.}
\label{fig::surrogate-pic}
\end{figure}

We give a schematic diagram in
Figure \ref{fig::surrogate-pic}.
The gray area represents
$\cS(\epsilon_2,\epsilon_\infty)$.
These are the functions that preclude adaptivity.
Being close to $\cF$ in $L^2$ makes them hard to detect
but being far from $\cF$ in $L^\infty$ makes them hard to cover.
To achieve adaptivity we must settle for sometimes
covering $\Pi_\cF f$.

\subsubsection{Lower Bounds}

We begin with two lemmas.
The first controls the minimum width of a band and
the second controls the maximum.
The second is of more interest for our purposes; the first lemma
is included for completeness.
For any $1 \le p \le \infty$, 
$\epsilon>0$, and $A\subset \mathbb{R}^n$ define
\begin{equation}
M_p(\epsilon,A) = \sup\{ d_{\srm TV}(P_f,P_g):\ f,g\in A, ||f-g||_p \le \epsilon \}
\end{equation}
and
\begin{equation}
m_\infty(\epsilon,A_0,A_1)=
\inf \{d_{\rm TV}(P_f,P_g):\ f\in A_0,g\in A_1, \norm{f-g}_\infty \ge \epsilon\}.
\end{equation}

\begin{lemma}\label{lemma::mod}
Suppose that  
$\inf_{f\in A}\P_f{L \le f \le U} \ge 1-\alpha$.
Let $1 \le p \le \infty$ and $\epsilon > 0$.
For $f\in A$, define 
$$\epsilon(f,q) = \sup\{ \norm{f-h}_q:\ h\in A, \norm{f-h}_p\le \epsilon\},$$
where $1 \le q \le \infty$.
Then, for any $A_0\subset A$,
\begin{equation} \label{eq::eps-}
\inf_{f\in A_0} \P_f{W > \epsilon(f,\infty)} \ge  1- 2\alpha - \sup_{f\in A_0} M_p(\epsilon(f,p),A)
\end{equation}
where $W = ||U-L||_\infty$.
If every point in $A$ is contained in a subset of $A$ of $\ell^p$-diameter $\epsilon$,
then $\epsilon(f,p) \equiv \epsilon$, and
\begin{equation}
\inf_{f\in A_0} \P_f{W > \epsilon} \ge  1- 2\alpha - M_p(\epsilon,A).
\end{equation}
\end{lemma}

\begin{lemma}\label{lemma::m2}
Suppose that  
$\inf_{f\in A}\P_f{L \le f \le U} \ge 1-\alpha$.
Suppose that $A=A_0\cup A_1$ (not necessarily disjoint).
Let $\epsilon > 0$ be such that
for each $f\in A_0$ there exists $g\in A_1$ for which
$\norm{f-g}_\infty = \epsilon$.
Then,
\begin{equation} 
\sup_{f\in A_0} \P_f{W > \epsilon} \ge  1- 2\alpha -  m_\infty(\epsilon,A_0,A_1)
\end{equation}
where $W = ||U-L||_\infty$.
\end{lemma}

Now we establish the target rate, the smallest width
of a band if we knew a priori that
$f\in \cF$.
Define
\begin{equation}\label{eq::define-target}
w_\cF \equiv w_\cF(\alpha,\gamma,\sigma) = \Omega_\cF \, \sigma\,\tau^{-1}(1-2\alpha - \gamma).
\end{equation}

\begin{theorem}\label{thm::target-rate}
Suppose that
\begin{equation}
\inf_{f\in \cF}\P_f{L \le f\le U} \ge 1-\alpha.
\end{equation}
If $\inf_{f\in\cF}\P_f{W\le w}\ge 1-\gamma$
then $w \ge w_{\cF}$.

A band that achieves this width, up to logarithmic factors, 
is $(L,U) = \hat{f} \pm c$
where $\hat{f} = \Pi Y$
and $c = \sigma (\Pi \Pi^T)_{ii} z_{\alpha/2n}$.
\end{theorem}

\begin{remark}
Using an argument similar to that in Theorem \ref{thm::basic},
it is possible to improve this lower bound by an additional
$\sqrt{\log d}$ factor, but this is inconsequential to the
rest of the paper.
\end{remark}

Next, we give the main result for this case.
\begin{eqnarray}
v_0(\epsilon_2,\epsilon_\infty, n, \alpha, \gamma, \sigma) &=& \min\Bigl\{\sqrt{n}\epsilon_2,\, \epsilon_\infty,\, \sigma\tau^{-1}(1-2\alpha-\gamma)\Bigr\}, \label{eq::v0const}\\
v_1(\epsilon_2, n, d, \alpha, \gamma, \sigma)
   &=&  \left\{
    \begin{array}{ll}
     0 & {\rm if\ } \epsilon_2 \ge 2v_2(n,d,\alpha,\gamma) \label{eq::v1const}\\
     v_2(n,d,\alpha,\gamma)& {\rm if\ }\epsilon_2 < 2v_2(n,d,\alpha,\gamma),
    \end{array}
    \right. \\
v_2(n,d,\alpha,\gamma) &=& \kappa(\alpha,\gamma)(n-d)^{1/4}n^{-1/2} \\
\noalign{\noindent and define}
v(\epsilon_2,\epsilon_\infty, n, d, \alpha, \gamma, \sigma) &=& \max( v_0, v_1 ). %\Biggl\{ v_0(\epsilon_2,\epsilon_\infty, n, \alpha, \gamma, \sigma), v_1(\epsilon_2, n, d, \alpha, \gamma, \sigma) \Biggr\}.
\end{eqnarray}

\smallskip

\begin{theorem}[Lower Bound for Surrogate Confidence Band Width] \label{thm::main} 
%\newline
\hfil\break
Fix $0< \alpha < 1$ and
$0 < \gamma < 1-2\alpha$.
Suppose that for bands $B = (L,U)$
\begin{equation} \label{eq::validband::singlespace}
\inf_{f\in\R^n} \P_f{ F^*(f) \cap B \ne \emptyset }\ge 1-\alpha. 
\end{equation}
Then,
\begin{equation}\label{eq::Wgamma}
\inf_{f\in\cF} \P_f{W \le w} \ge 1 - \gamma.
\end{equation}
implies
\begin{equation}\label{eq::the-cool-lower-bound}
w \ge \underline{w}(\cF,\epsilon_2,\epsilon_\infty,n,d,\alpha,\gamma,\sigma) \equiv \max\Bigr\{w_\cF(\alpha,\gamma,\sigma),\, v(\epsilon_2,\epsilon_\infty, n, d, \alpha, \gamma, \sigma)\Bigr\}
\end{equation}
\end{theorem}

The inequality (\ref{eq::validband::singlespace}) ensures that $B$ is a valid surrogate confidence band: for every function,
either the function or its surrogate is covered with at least the target probability.
The result gives a probabilistic lower bound on the width of the band that is at least as big as the best a priori width for the subspace.
As we will see, with proper choice of $\epsilon_2$ and $\epsilon_\infty$, the $v$ term can be made small, giving the subspace width $w_\cF$
for the lower bound.

Next, we address the question of optimality. Consider, for example, the trivial surrogate that maps all functions to 0.
We can cover the surrogate using 0 width bands with probability 1, but this would not be too interesting.
There is a tradeoff between the width of the bands on low dimensional subspaces and the volume of the spoiler
set, the functions that are surrogated.
We characterize optimality here as minimizing the volume of the spoiler set $\cS(\epsilon_2,\epsilon_\infty)$
while still attaining the target width with high probability when $f$ truly lies in the subspace.
In this sense, the surrogate defined above is optimal.

\begin{theorem}[Optimality]
\label{thm::optimal}
Let $\underline{w}$ denote the 
right hand side of inequality
(\ref{eq::the-cool-lower-bound}).
Then $\underline{w} \ge w_\cF$,
where $w_\cF$ is defined in (\ref{eq::define-target}).
Setting
$$
\epsilon_2 = 2\kappa(\alpha,\gamma)(n-d)^{1/4}n^{-1/2},\ \ \ 
\epsilon_\infty = w_\cF
$$
minimizes ${\rm Volume}(\cS(\epsilon_2,\epsilon_\infty))$ subject to
achieving the lower bound on $\underline{w}$.
\end{theorem}

\subsubsection{Achievability}

Having established a lower bound, we need to show that the lower bound is sharp.
We do this by constructing a finite-sample procedure that achieves the bound within
a factor of 2.
Let
$F_{a,d}$ denote the {\sc cdf} of a $\chi^2$ 
random variable with $d$ degrees of
freedom and noncentrality parameter $a$
and let
$\chi^2_{\alpha,d}= F^{-1}_{0,d}(1-\alpha)$.
Let $T = ||Y-\Pi Y||^2$ and
define
\begin{equation}\label{eq::these-are-the-bands}
B=(L,U) = \hat{f}\pm c \sigma
\end{equation}
where
\begin{equation}
\hat{f} =
\left\{
\begin{array}{ll}
Y      & {\rm if\ } T > \chi^2_{\gamma,n-d}\\
\Pi Y  & {\rm if\ } T \le \chi^2_{\gamma,n-d}
\end{array}
\right.
\end{equation}
and
\begin{equation}
c = z_{\alpha/2n} \times 
\left\{
\begin{array}{ll}
\omega_\cF + \epsilon_\infty & {\rm if}\  T \le \chi^2_{\gamma,n-d}\\
1                    & {\rm if}\  T > \chi^2_{\gamma,n-d}.
\end{array}
\right.
\end{equation}

\begin{theorem}\label{thm::achieve}
If
\begin{equation}
\gamma \ge 1- F_{0,n-d}(F^{-1}_{n\epsilon_2^2,n-d}(\alpha/2))
\end{equation}
then
\begin{equation}
\inf_{f\in \mathbb{R}^n} \P_f{F^\star(f) \cap B \ne \emptyset} \ge 1-\alpha
\end{equation}
and
\begin{equation}\label{eq::it-works}
\inf_{f\in \cF} \P_f{W \le  w_\cF + \epsilon_\infty} \ge 1-\gamma.
\end{equation}
If $\epsilon_2 \ge E(n-d,\alpha/2,\gamma) (n-d)^{1/4} n^{-1/2}$,
where $E(m,\alpha,\gamma)$ is defined in (\ref{eq::const::E}),
then
\begin{equation}\label{eq::it-works2}
\inf_{f\in \cF} \P_f{W \le  2 \underline{w}(\cF,\epsilon_2,\epsilon_\infty,\alpha,\gamma,n,d)}
\ge 1-\gamma.
\end{equation}
where
$\underline{w}(\cF,\epsilon_2,\epsilon_\infty,\alpha,\gamma,n,d)$
is defined (\ref{eq::the-cool-lower-bound}).
Hence,
the procedure adapts to within a logarithmic factor
of the lower bound $\underline{w}$
given in Theorem \ref{thm::main}.
\end{theorem}

\begin{corollary}
Setting
$$
\epsilon_2 = E(n-d,\alpha/2,\gamma)(n-d)^{1/4}n^{-1/2},\ \ \ 
\epsilon_\infty = w_\cF
$$
in the above procedure,
minimizes ${\rm Volume}(\cS(\epsilon_2,\epsilon_\infty))$ subject to
satisfying (\ref{eq::it-works2}).
\end{corollary}

\begin{remark}\label{remark::sigma}
The results can be extended to unknown $\sigma$
by replacing $\sigma$ with a nonparametric estimate $\hat\sigma$.
However, the results are then asymptotic rather than finite sample.
Moreover, a minimal amount of smoothness is required
to ensure that $\hat\sigma$ consistently estimates $\sigma$;
see Genovese and Wasserman (2005).
So as not to detract from our main points,
we continue to take $\sigma$ known.
\end{remark}

\subsubsection{Remarks on Estimation and the Modulus of Continuity}

It is interesting to note that the bands 
defined above cover the true $f$
over a set $V$ that is larger than $\cF$.
In this section we take a brief look at the properties of $V$.

Define 
\begin{equation}
C(\alpha,a,b) = \sup_{u > 0}\,(a u + b)\left(1 - \alpha - \frac14 + \half \Phi(-u/2)\right),
\end{equation}
and let $C(\alpha) \equiv C(\alpha,1,0)$.
Let $\cF^\perp$ be
the orthogonal complement of $\cF$. Let $B^\perp_k(0,\epsilon)$
be a $\ell^k$-ball around 0 in $\cF^\perp$ ($k = 2,\infty$). 
For $f\in\mathbb{R}^n$, let $B^\perp_k(f,\epsilon) = f + B^\perp_k(0,\epsilon)$.
Define 
\begin{equation}\label{eq::V}
V \equiv V(\cF,\epsilon_2,\epsilon_\infty) = 
\bigcup_{f\in\cF} \Biggl(B^\perp_2(f,\epsilon_2) \cap B^\perp_\infty(f,\epsilon_\infty)\Biggr).
\end{equation}

\begin{lemma}
Let $B=(L,U)$ be defined as in (\ref{eq::these-are-the-bands}).
Then
\begin{equation}
\inf_{f\in V}\P_f{L \leq f \leq U}\geq 1- \alpha.
\end{equation}
\end{lemma}

Let $T f = f_1$. The next lemma gives the modulus of continuity 
(Donoho and Liu 1991)
of $T$ over $V$
which measures the difficulty of estimation over $V$.
The modulus of continuity of $T$ over a set $\cA$ is
\begin{equation}
\omega(u,\cA) = \sup\{ |T f - T g| :\; \norm{f-g}_2 \le u; f, g\in \cA\}.
\end{equation}
%\begin{eqnarray}
%\omega(u,\cA) &=& \sup\{ |T f - T g| :\; \norm{f-g}_2 \le u; f, g\in \cA\}\\
%&=& \sup\{ |T f - T g| :\; \norm{f-g}_\natural \le u\sqrt{n}; f, g\in \cA\}
%\end{eqnarray}
%where $||x||_\natural = \sqrt{\sum_i x_i^2}$.
Donoho and Liu showed that the difficulty of estimation over $\cA$
is often characterized by $\omega(1/\sqrt{n},\cA)$
in the sense that this quantity defines a lower bound
on estimation rates.

\begin{lemma}[Modulus of Continuity]
\label{lemma::Modulus-of-Continuity}
We have
\begin{equation}
\omega(u,V) = 
\left( u \Omega\sqrt{n} \sqrt{\frac{\Omega^2}{1 + \Omega^2}} + 
\min\left( \frac{u\sqrt{n}}{\sqrt{1+\Omega^2}}, 
\epsilon_2\wedge (\epsilon_\infty/\sqrt{n})\right)\right).
\end{equation}
\end{lemma}

Note that when $\epsilon_2 = \epsilon_\infty =0$ and
$\Omega \sim \sqrt{d/n}$, we have
$\omega(1/\sqrt{n},\cA)\sim \sqrt{d/n}$
as expected.
However, when
$\epsilon \equiv \epsilon_2 = \epsilon_\infty/\sqrt{n}$ is large
we will have that
$\omega(1/\sqrt{n},\cA)\sim \sqrt{d/n} + \epsilon/\sqrt{1+d^2/n}$.
The extra term
$\epsilon/\sqrt{1+d^2/n}$ reflects the ``ball-like'' behavior
of $V$ in addition to the subspace-like behavior of $V$.
The bands need to cover over this extra set to maintain valid coverage
and this leads to larger lower bounds than just covering over $\cF$.

%
%The following lemma follows directly from Theorem 1 of Low (1997).
%
%\begin{lemma}\label{lemma::from-low}
%Let $X$ be multivariate Normal with mean $\mu\in\mathbb{R}^n$
%and covariance $\sigma^2 I$ where
%$I$ is the $n\times n$ identity matrix.
%Let $C(X)$ be a $1-\alpha$ confidence interval for
%$T\mu$ 
%where $T$ is a linear functional.
%Let $W$ be the width of $C(X)$.
%If $\cF$ is convex
%then
%\begin{equation}
%\sup_{f\in\cF}\E_f (W) \ge 
%\sup_{\xi > 0}\,\left(1 - \alpha - \frac\xi{4}\right)
%\omega\left( -2 \Phi^{-1}\left(\frac{1-\xi}2\right), \cF \right).
%\end{equation}
%\end{lemma}
%
%
%\begin{corollary} If $\cF$ is a subspace and $(L,U)$ a valid $1-\alpha$ confidence
%band over $\cF$, then
%\begin{equation}
%\sup_{f\in\cF}\E_f W \ge w(\cF,\alpha).
%\end{equation}
%where $w(\cF,\alpha) = C(\alpha) \Omega_{\cF}. $
%\end{corollary}
%
%
%\begin{theorem}
%If $(L,U)$ a valid $1-\alpha$ confidence band over that $V$
%and $\epsilon = \max(\epsilon_{2,\natural},\epsilon_\infty)$, then
%\begin{equation}
%\sup_{f\in V}\E_f (W) \ge 
%\max\left( C(\alpha) \sqrt{1 + \Omega^2_{\cF}},\, 
%C\left(\alpha, \Omega_{\cF} \sqrt{\frac{\Omega^2_{\cF}}{1 + \Omega^2_{\cF}}}, \epsilon\right) \right).
%\end{equation}
%\end{theorem}
%
%
%
%ATTN: NEED SOME COMMENTS HERE
%
%

\subsection{Nested Subspaces} \label{sec::nestedspaces}

Now suppose that we have nested subspaces
$\cF_1 \subset \cdots \subset \cF_m \subset \cF_{m+1}\equiv \mathbb{R}^n$.
Let $\Pi_j$ denote the projector onto $\cF_j$.
We define the surrogate as follows.

\begin{definition}
For given 
$\epsilon_2 = (\epsilon_{2,1}, \ldots, \epsilon_{2,m})$
and
$\epsilon_\infty = (\epsilon_{\infty,1}, \ldots, \epsilon_{\infty,m})$
define
\begin{equation}
\cJ(f) = \{1\le j\le m:\ 
||f-\Pi_j f||_2 \le \epsilon_{2,j}\ {\rm and}\ ||f-\Pi_j f||_\infty > \epsilon_{\infty,j}
\Biggr\}.
\end{equation}
Then define the surrogate set
\begin{equation}
F^\star(f) = \{ \Pi_j f:\ j\in \cJ(f)\} \cup \{f\}.
\end{equation}
\end{definition}

\begin{definition}
We say that $B=\{ g:\ L\le g \le U\}\equiv (L,U)$ has coverage $1-\alpha$ if
\begin{equation}
\inf_{f\in\mathbb{R}^n}\P_f{ F^\star\cap B \neq \emptyset}\ge 1-\alpha.
\end{equation}
\end{definition}

\subsubsection{Lower Bounds}

\begin{theorem}[Lower Bound for Surrogate Confidence Band Width] \label{thm::main::nested} 
%\newline
\hfil\break
Fix $0< \alpha < 1$ and
$0 < \gamma < 1-2\alpha$.
Suppose that for bands $B = (L,U)$
\begin{equation} \label{eq::validband::nestedspaces}
\inf_{f\in\R^n} \P_f{ F^*(f) \cap B \ne \emptyset }\ge 1-\alpha. 
\end{equation}
Then 
\begin{equation}\label{eq::Wgamma::nested}
\inf_{f\in\cF_j} \P_f{W \le w} \ge 1 - \gamma.
\end{equation}
implies
\begin{equation}\label{eq::the-cool-lower-bound::nested}
w \ge \underline{w}(\cF_j,\epsilon_{2,j},\epsilon_{\infty,j},n,d_j,\alpha,\gamma,\sigma),
\end{equation}
where $\underline{w}$ is given in Theorem \ref{thm::main}.
\end{theorem}

\begin{theorem}[Optimality]
\label{thm::optimal::nested}
Let $\underline{w}$ denote the 
right hand side of inequality
(\ref{eq::the-cool-lower-bound::nested}).
Then $\underline{w} \ge w_\cF$,
where $w_{\cF_j}$ is defined in (\ref{eq::define-target}).
Setting
$$
\epsilon_{2_j} = 2\kappa(\alpha,\gamma)(n-d_j)^{1/4} n^{-1/2},\ \ \ 
\epsilon_{\infty,j} = w_{\cF_j}
$$
minimizes the volume of the set
\begin{equation} \label{eq::minvolume::nested}
\Set{ f:\ \norm{f - \Pi_j f} \le \epsilon_{2,j} \ {\rm and}\ \norm{f - \Pi_j f}_\infty > \epsilon_{2,\infty}}
\end{equation}
subject to achieving the lower bound on $\underline{w}$.
\end{theorem}

\subsubsection{Achievability}

Define
$T_j= ||Y-\Pi_j Y||^2$ and $\hat{f}=\Pi_{\hat J}Y$, where
\begin{equation}
\hat{J} = \min\{1 \le j \le m:\ T_j \le \chi^2_{\gamma,n-d_j}\},
\end{equation}
where $\hat J = m + 1$ if the set is empty, 
and define
\begin{equation}
c_j = z_{\alpha_j/2 n} \times 
\left\{
\begin{array}{ll}
\omega_{\cF_j}(\alpha_j) + \epsilon_{\infty,j} & {\rm if}\  1 \le j \le m\\
 1                                             & {\rm if}\ j = m + 1.
\end{array}
\right.
\end{equation}
Finally, let $B=(L,U) = \hat{f}\pm c_{\hat{J}} \sigma$
where $\sum_j \alpha_j \le \alpha$.

\begin{theorem}\label{thm::achieve-multi}
If, 
\begin{equation}\label{eq::gam-alph}
\gamma \ge 1- \min_j F_{0,n-d_j}(F^{-1}_{n\epsilon_{2,j}^2,n-d_j}(\alpha_j))
\end{equation}
then
\begin{equation}
\inf_{f\in\mathbb{R}^n}\P_f{ F^\star\cap B \neq \emptyset}\ge 1-\alpha.
\end{equation}
Let $w_j = w_{\cF_j}(\alpha_j) + \epsilon_{\infty,j}$.
If $w_1 \le \cdots \le w_{m+1}$ then
\begin{equation}\label{eq::the-gam-result}
\inf_{f\in \cF_j} \P_f{W \le  w_j} \ge 1-\gamma.
\end{equation}
If in addition 
$\epsilon_{2,j} \ge E(n-d_j,\alpha_j,\gamma)(n-d_j)^{1/4} n^{-1/2}$
and
$\epsilon_{\infty,j} \le w_{\cF_j}$
then
\begin{equation}\label{eq::opt-this}
\inf_{f\in \cF_j} 
\P_f{W \le  2 \underline{w}(\epsilon_{2,j},\epsilon_{\infty,j},\alpha_j,\gamma,n,d_j)}
\ge 1-\gamma
\end{equation}
where
$\underline{w}(\epsilon_{2,j},\epsilon_{\infty,j},\alpha_j,\gamma,n,d_j)$
is defined (\ref{eq::the-cool-lower-bound}).
Hence,
the procedure adapts to within a logarithmic factor
of the lower bound $\underline{w}$
given in Theorem \ref{thm::main}.
\end{theorem}

\begin{corollary}\label{cor::optimality}
Suppose $\alpha_1 = \cdots = \alpha_{m+1} = \alpha/(m+1)$.
Then
$w_1 \le \cdots \le w_{m+1}$
so (\ref{eq::the-gam-result}) holds.
Moreover,
setting
\begin{eqnarray}
\epsilon_{2,j} &=& E(n-d_j,\alpha_j,\gamma)(n-d_j)^{1/4}n^{-1/2}\\
\noalign{\noindent and}
\epsilon_{\infty,j} &=& w_{\cF_j}
\end{eqnarray}
in the above procedure,
minimizes the volume of the set (\ref{eq::minvolume::nested})
satisfying (\ref{eq::the-cool-lower-bound::nested}).
\end{corollary}

\begin{example}
Suppose that
$x_i=i/n$ and let
$B_1=[0,1/d], B_2 = (1/d,2/d]$, $\ldots$, $B_d = ((d-1)/d,1]$.
Write
$f =( f(x_i):\ i=1, \ldots, n)$
and let $\cF$
denote the subspace of vectors $f$
that are constant over each $B_j$.
Then $\Omega_{\cF} = \sqrt{d/n}$.
The above procedure
then produces a band with width no more that
$O(\sqrt{d/n})$ with probability at least $1-\gamma$.
\end{example}

\section{Proofs}
\label{sec::proofs}

In this section, we prove the main results.
We omit proofs for a few of the simpler lemmas.
Throughout this section, we write
$x_n = O^*(b_n)$ to mean that
$x_n = O(c_n b_n)$ 
where $c_n$ increases at most logarithmically with $n$.

The following lemma is essentially from Section 3.3 of Ingster and Suslina (2003).

\begin{lemma}\label{lemma::mix}
Let $M$ be a probability measure
on $\mathbb{R}^n$ and let
$$
Q(\cdot) = \int P_f(\cdot) dM(f)
$$
where
$P_f(\cdot)$ denotes the measure for a multivariate
Normal with mean $f = (f_1, \ldots, f_n)$
and covariance $\sigma^2 I$.
Then
\begin{equation}
L_1(Q,P_g) \le 
\sqrt{\int \int \exp\left\{ \frac{n \langle f-g, \nu-g\rangle}{\sigma^2}\right\}
dM(f)dM(\nu) - 1}.
\end{equation}
In particular, if $Q$ is uniform on a finite set $\Omega$, then
\begin{equation}
L_1(Q,P_g) \le \sqrt{
\left( \frac{1}{|\Omega|}\right)^2
\sum_{f,\nu\in\Omega} \exp\left\{ \frac{n \langle f-g, \nu-g\rangle}{\sigma^2}\right\}-1}.
\end{equation}
\end{lemma}

\begin{proof}%[ of Lemma \ref{lemma::mix}]
Let $p_f$ denote the density
of a multivariate Normal with mean
$f$ and covariance $\sigma^2 I$
where $I$ is the identity matrix.
Let $q$ be the density of $Q$:
$$
q(y) = \int p_f(y) dM(f).
$$
Then,
\begin{eqnarray}
\nonumber
\int | p_g(x) - q(x) | dx &=&
  \int \frac{| p_g(x) - q(x) |}{\sqrt{p_g(x)}} \sqrt{p_g(x)} dx \\
\label{eq::ineq}
  &\le&  \sqrt{\int  \frac{( p_g(x) - q(x) )^2}{p_g(x)}  dx} =
       \sqrt{\int \frac{q^2(x)}{p_g(x)} dx -1 }.
\end{eqnarray}
Now,
\begin{eqnarray*}
\int \frac{q^2(x)}{p_g(x)} dx  &=& 
\int \left(\frac{q(x)}{p_g(x)}\right)^2 p_g(x) dx = \E_g  \left(\frac{q(x)}{p_g(x)}\right)^2 \\
&& \hskip -2cm = \int \int  \E_g\left(\frac{p_f(x)p_\nu(x)}{p_g^2(x)}\right)dM(f)dM(\nu)\\
&& \hskip -2cm=\int\int \exp\left\{ - \frac{n}{2\sigma^2} ( ||f-g||^2 + ||\nu-g||^2 )\right\}
     \E_g \left(\exp\left\{ \epsilon^T (f+\nu-2g)/\sigma^2\right\}\right)dM(f)dM(\nu)\\
&& \hskip -2cm=\int\int \exp\left\{ - \frac{n}{2\sigma^2} ( ||f-g||^2 + ||\nu-g||^2 )\right\}
     \exp\left\{ \sum_{i=1}^n (f_i - g_i + \nu_i - g_i)^2 /(2\sigma^2) \right\}dM(f)dM(\nu)\\
&& \hskip -2cm=\int\int \exp\left\{ \frac{n\langle f-g, \nu-g\rangle}{\sigma^2}\right\}dM(f)dM(\nu)
\end{eqnarray*}
and the result follows from (\ref{eq::ineq}).
\end{proof}

\begin{proof}[Proof of Theorem \ref{thm::basic}]
Let $N = |\Omega|$ and
let $b^2 = n\max_{f\in\Omega}||f-g||^2$.
Let $p_f$ denote the density
of a multivariate Normal with mean
$f$ and covariance $\sigma^2 I$
where $I$ is the identity matrix.
Define the mixture
$$
q(y) = \frac{1}{N}\sum_{f\in\Omega}p_f(y).
$$
By Lemma \ref{lemma::mix},
\begin{eqnarray*}
\int | p_g(x) - q(x) | dx    &\le&
\sqrt{ \left(\frac{1}{N}\right)^2
\sum_{f,\nu\in\Omega} \exp\left\{ \frac{n\langle f-g, \nu-g\rangle}{\sigma^2}\right\} -1}\\
&=&
\sqrt{ \left(\frac{1}{N}\right)^2 \Biggl[ N e^{b^2/\sigma^2} + N(N-1)\Biggr]-1}\\
  &\le&  \sqrt{ e^{b^2/\sigma^2}/N} = \epsilon.
\end{eqnarray*}

Define two events, 
$A=\{ \ell \le g \le u\}$ and
$B = \{ \ell \le f \le u, \ {\rm for\ some \ }f\in \Omega\}$.
Then,
$A \cap B \subset \{ w_n \ge a\}$
where
$$
a= \min_{f\in\Omega}||g-f||_\infty.
$$
Since
$\P_f{\ell \le f\le u}\ge 1-\alpha$ for all $f$,
it follows that
$\P_f{B} \ge 1-\alpha$ for all $f\in\Omega$.
Hence,
$Q(B) \ge 1-\alpha$.
So,
\begin{eqnarray*}
\P_g{w_n \ge a}   &\ge& \P_g{A\cap B} \ge  Q(A \cap B) - \epsilon  =  Q(A) + Q(B) - Q(A\cup B)- \epsilon\\
&\ge & Q(A) + Q(B) - 1- \epsilon \ge  Q(A) + (1-\alpha) - 1- \epsilon \ge  \P_g{A} + (1-\alpha) - 1- 2\epsilon\\
  &\ge& (1-\alpha) + (1-\alpha) - 1- 2\epsilon = 1-2\alpha - 2\epsilon.
\end{eqnarray*}
So,
$\E_g(w_n) \ge (1-2\alpha - 2\epsilon) a$.
\end{proof}

\begin{proof}[Proof of Theorem \ref{thm::finite-sample-lower-bound}]
Let $g\in\mathbb{R}^n$ be arbitrary, let
$$
a_n = \sigma \sqrt{\log( n \epsilon^2)}
$$
and define
$$
\Omega = \Biggl\{ g+(a_n,0,\ldots,0), \ g+(0,a_n,\ldots,0), \ \ldots, \ g+ (0,0,\ldots,a_n) \Biggr\}.
$$
Then
the conditions of Theorem \ref{thm::basic} are satisfied
with $N=n$, and
hence
\begin{equation}
\E_g(W)\ge  (1-2\alpha - 2\epsilon)\min_{f\in\Omega}||g-f||_\infty =
(1-2\alpha - 2\epsilon) a_n.
\end{equation}
This is true for each $g$ and hence
(\ref{eq::rnbound}) follows.
The last statement of the theorem follows from
standard Gaussian tail inequalities.
\end{proof}

\begin{proof}[Proof of Theorem \ref{thm::finite-sample-lip}]
We construct the appropriate set $\Omega$
and apply Theorem \ref{thm::basic}.
For simplicity, we build $\Omega$ around
$g=(0,\ldots, 0)$, the extension to arbitrary $g$ being straightforward.
Set $a = a_n$ from the statement of the theorem, and 
define
$$
F(x) = \left\{
\begin{array}{ll}
Lx       & 0 \le x \le a/L \\
2a -Lx   & a/L \le x \le 2 a /L.
\end{array}
\right.
$$
Note that 
$F\in {\cal F}(L)$
and that $F$ minimizes $||F||_2$ among all $F\in  {\cal F}(L)$
with $||F||_\infty =a$.
For simplicity, assume that
$2aN/L = 1$ for some integer $N$.
Define
$F_1(\cdot) = F(\cdot)$,
$F_2(\cdot) = F(\cdot-\delta)$,\ldots,
and $F_N(\cdot) = F(\cdot-N\delta)$.
Let
$\Omega(a) = \{ f_1, \ldots, f_N\}$
where
$f_j = (F_j(x_1), \ldots, F_j(x_n))$.
Now
$$
n ||f_j||^2 \le  \frac{2 n a^3}{3L}
$$
and so
$$
\frac{e^{n ||f_j||^2/\sigma^2}}{N} \le \epsilon^2.
$$
Now apply Theorem \ref{thm::basic}. 

To prove the last statement,
we note that it is well known that if $\hat{F}$ 
is a kernel estimator with triangular kernel and bandwidth
$h=O(n^{-1/3})$ then
$$
\sup_{f\in \Theta} E_F (||\hat{F} - F||_\infty) \le 
C \left(\frac{ \log n}{n}\right)^{1/3} \equiv C_n 
$$
for some $C>0$.
Then $B=(\hat{F}-\frac{C_n}{\alpha},\hat{F}+\frac{C_n}{\alpha})$
(restricted to $x_i=i/n$) is valid by Markov's inequality
and has the rate $a_n$.
\end{proof}

\begin{proof}[Proof outline of Theorem \ref{thm::sobolev}]
We will use the fact that an appropriately chosen wavelet basis
forms a basis for $\cF$.
Let
$$
J_n \sim \log_2 \left(  \frac{n^{1/(2p+1)}}{\log n} \right),
$$
$$
b_n = \frac{\sigma}{\sqrt{n}}\sqrt{\log (2^{J_n}\epsilon^2)}
$$
and
$$
F(x) = b_n 2^{J_n/2}\psi(2^{J_n}x)
$$
where
$\psi$ is a compactly supported mother wavelet.
Then
$F^{(p)} = b_n 2^{J_n/2}2^{p J_n}\psi^{(p)}(2^{J_n}x)$ so that
$\int (F^{(p)})^2 < c^2$ for all large $n$
so that $F\in\cF$.

Let $f=(F(x_i),\ldots, F(x_n))$.
Then,
$$
||f||_\infty = b_n 2^{J_n/2} = O^*(n^{-p/(2p+1)})
$$
and
$\sqrt{n}||f||_2 \sim \sqrt{n} b_n$.
Let
$f_k = (F(x_1 - k\Delta),\ldots, F(x_n - k\Delta))^T$
where $\Delta$ is just large enough so that
the $F_k$'s are orthogonal.
Hence,
$\Delta \approx 1/N$ where
$N \sim 2^{J_n}$.
Finally, set
$\Omega = \{ f_1, \ldots, f_N\}$.
Then,
$$
\frac{e^{n ||f||^2/\sigma^2}}{N} = e^{n b_n^2/\sigma^2}{2^{J_n}} \le \epsilon^2
$$
for each $f\in\Omega$.
The lower bound follows from Theorem \ref{thm::basic}.

A fixed-width procedure that achieves the bound is
$$
\ell_i = \hat{f}_i -  c_n z_{\alpha/n},\ u_i = \hat{f}_i +  c_n z_{\alpha/n}.
$$
where
$\hat{f}_i = \hat{F}(x_i)$,
$$
\hat{F}(x) = \sum_j \hat\alpha_j \phi_j(x) + \sum_{j=1}^J \sum_k \hat\beta_{jk}\psi_{jk}(x),
$$
$\hat\alpha_j = n^{-1}\sum_i Y_i\phi_j(x_i)$,
$\hat\beta_{jk}= n^{-1}\sum_i Y_i\psi_{jk}(x_i)$
and
$c_n = \sqrt{\max_x {\rm Var}(\hat{F}(x))}$.
\end{proof}

\begin{proof}[Proof outline of Theorem \ref{thm::besov}]
Again, we use the fact that an appropriately chosen wavelet basis
forms a basis for $\cF$.
Let
$$
J_n \sim \frac{ \log_2 \frac{c \sqrt{n}}{\sigma \sqrt{\log 2^J \epsilon^2}}}
{\xi + \frac{1}{2} - \frac{1}{p}}.
$$
Let 
$$
a_n = \frac{\sigma}{\sqrt{n}} \sqrt{ \log 2^J \epsilon^2}
$$
and define
$F(x) = a_n 2^{J/2} \psi(x)$,
where
$\psi$ is a compactly supported mother wavelet.
Then,
$||f|| =a_n$,
$||f||_\infty =a_n 2^{J/2}$,
and
$||F||_{p,q}^\xi \le c-\delta$ for all large $n$.
Take $\Omega$ around $g$ to be non-overlapping translations of $F$ added to $g$.
Then $N \sim 2^J$
and conditions of Theorem \ref{thm::basic} hold.
Moreover,
$$
a_n = O^*( n^{-1/(1/p - \xi - 1/2)}).
$$
The bound is achieved
by Markov applied to the soft-thresholded wavelet estimator 
with universal thersholding.
\end{proof}

\begin{proof}[Proof of Lemma \ref{lemma::universal::Q}]
$Q$ is the solution, with respect to $c$, to
$\xi = 1- F_{0,m}(r(c))$
where
the function
$r(c)=F^{-1}_{c\sqrt{m},m}(\beta))$
is monotonically increasing in $c$.
Also, $F_{0,m}(r(0))=\beta$ and 
$F_{0,m}(r(\infty))=1$
so a solution exists since
$0 < \beta < 1-\xi < 1$.
Now we bound $Q$ from above.

To upper bound $Q$ it suffices to find $c$ such that
\begin{equation}
F^{-1}_{c\sqrt{m},m}(\beta) \ge F^{-1}_{0,m}(1-\xi).
\end{equation}
\relax From Birg\'e (2001) we have
\begin{eqnarray}
F^{-1}_{z,d}(u)   &\le& z+d+2\sqrt{(2z+d)\log(1/(1-u))} + 2\log (1/(1-u))\\
F^{-1}_{z,d}(u)   &\ge& z+d-2\sqrt{(2z+d)\log(1/u)}.
\end{eqnarray}
Hence,
\begin{eqnarray}
F^{-1}_{c\sqrt{m},m}(\beta) &\ge& m + c\sqrt{m} - 2 \sqrt{(2 c \sqrt{m} + m) \log\frac1\beta} \\
F_{0,m}^{-1}(1-\gamma) &\le& m + 2 \sqrt{m \log\frac1\gamma} + 2\log\frac1\gamma.
\end{eqnarray}
It suffices to find $c$ that satisfies
\begin{equation}
m + c\sqrt{m} - 2 \sqrt{(2 c \sqrt{m} + m) \log\frac1\beta} 
  \ge m + 2 \sqrt{m \log\frac1\gamma} + 2\log\frac1\gamma, \\
\end{equation}
or equivalently,
\begin{equation}
c \ge 2\sqrt{\left(\frac{c}{\sqrt{m}} + 1\right)\log\frac1\beta} + 2\left(\sqrt{\log\frac1\gamma} + \log\frac1\gamma\right).
\end{equation}
The right hand side of the last inequality is largest when $m = 1$, and
equality can be achieved when $m = 1$ at some $\Lambda(\beta,\xi)$
for any $\beta, \xi$ satisfying the stated conditions.
Equality can be achieved then for any $m$ at some $Q(m,\beta,\xi) \le \Lambda(\beta,\xi)$.
This proves the first claim.
The second claim follows immediately by inspection.
\end{proof}

\begin{proof}[Proof of Lemma \ref{lemma::min2inInfinity}]
Note that
\begin{eqnarray}
\min\Biggl\{ \norm{v}:\ v\in\cF, \norm{v}_\infty = 1 \Biggr\} &=&
\min_{v\in\cF}\frac{\norm{v}_{\phantom{\infty}}}{\norm{v}_\infty}\\
&=& \frac{1}{\max_{v\in\cF}\frac{\norm{v}_\infty}{\norm{v}_{\phantom{\infty}}}},\\
&=& \frac{1}{\max\Biggl\{ \norm{v}_\infty:\ v\in\cF, \norm{v} = 1 \Biggr\}}.
\end{eqnarray}
If $v$ solves one of these problems then
$\epsilon v$ solves the more general version in the statement of the lemma.
It now suffices to show just the second equality.

Now,
$\Omega_{\cF} = \max_i \Omega_i$
where
$$
\Omega_i = \frac{\langle e_i, \Pi_\cF  e_i\rangle}{\norm{e_i}\,\norm{\Pi_\cF e_i}} = 
\frac{\norm{\Pi_\cF    e_i}}{\norm{e_i}}.
$$
Maximizing $f_i=e_i^T f$
for $f\in\cF$ and $\norm{f}\le 1$ is equivalent to maximizing
$n \langle e_i, f\rangle = n \langle \Pi_\cF e_i, f\rangle$.  The
maximum subject to the constraint occurs at $f^\star =  \Pi e_i/\norm{\Pi e_i}$.
Hence, the maximum is
$e_i^T f^\star = (\Pi e_i)^T f^\star =
n\norm{\Pi e_i}^2/\norm{\Pi e_i} =
n\norm{\Pi e_i}^2/\norm{\Pi e_i} \frac{\norm{e_i}}{\norm{e_i}}=
\sqrt{n}\Omega_i$.
Maximizing over $i$ completes the proof.
\end{proof}

\begin{proof}[Proof of Lemma \ref{lemma::baraud-doover}]

We find a $P_0\in \cF_j$
and a measure $\mu$ supported on $A$ such that
$d_{\srm TV}(P_0, P_\mu) \le 2\delta$.
We then have, following Ingster (1993), 
\begin{eqnarray}
\beta 
    &\ge& \inf_{\phi_\xi\in\Phi_\xi} P_\mu\Event{ \phi_\xi = 0 } \\
    &\ge& 1 - \xi - \sup_{R:\ P_0(R) \le \xi} \left| P_0(R) - P_\mu(R) \right| \\
    &\ge& 1 - \xi - \sup_R \left| P_0(R) - P_\mu(R) \right| \\
    &=&  1 - \xi - \half d_{\srm TV}(P_0,P_\mu) \\
    &\ge& 1 - \xi - \delta.
\end{eqnarray}

Let $\psi_1, \psi_2, \ldots, \psi_{n}$ be an orthonormal basis for $\R^n$ such 
that $\psi_1,\ldots,\psi_d$ form an orthonormal basis for $\cF$.
Fix $\tau > 0$ small and let $\lambda^2 = n\epsilon^2/(n-d) + \tau^2/(n-d)$.
Define
\begin{equation}\label{eq::fE}
   f_E = \lambda \sum_{s=d+1}^{m} E_s \psi_s,
\end{equation}
where $(E_s: \ s=d+1,\ldots, n)$ are independent Rademacher random variables,
that is, $\P{E_s =1} = \P{E_s =-1} = 1/2$.
Now,
$\Pi_\cF f_E = 0$  and hence $||f_E - \Pi_\cF f_E||^2 = \lambda^2 > \epsilon^2$,
and hence $f_E\in A$ for each choice of the Rademachers.

Let $P_\mu = \E(P_E)$
where $P_E$ is the distribution under $f_E$ 
and the expectation is with repect to the Rademachers.
Choose $f_0 \in \cF$ and let $P_0$ be the corresponding
distribution.
As in Baraud, we use the bound
\begin{equation}
d_{\srm TV}(P_\mu,P_0) \le \sqrt{ \E_0 \left(\frac{dP_\mu}{d P_0}(Y)\right)^2 - 1 }.
\end{equation}
We take $f_0 = (0,\ldots, 0) \in \cF$ and so
\begin{eqnarray}
  \left(\frac{dP_\mu}{d P_0}(Y)\right)
     &=& \E_E \left(\exp\left\{-\frac{1}{2}\lambda^2(n-d) + 
                    \lambda\sum_{s=d+1}^{n} E_s \sum_i Y_i \psi_{si}\right\}\right)\\
     &=& e^{-\lambda^2/2}\,\prod_{s=d+1}^{n} \cosh(\lambda (Y\cdot \psi_s)).
\end{eqnarray}
Since
$\E_0 \cosh^2 (\lambda (Y\cdot\psi_j)) = e^{\lambda^2}\cosh(\lambda^2)$ 
and $\cosh(x) \le e^{x^2/2}$
we have
\begin{eqnarray}
  \E_0 \left(\frac{dP_\mu}{d P_0}(Y)\right)^2 
      &=&   \left(\cosh(\lambda^2)\right)^{n-d} \\
      &\le& e^{(n-d)\lambda^4/2} \\
      &=&   \exp\left(\frac{n^2}{2(n-d)} \epsilon^4 + \frac{\tau^4}{2(n-d)} + \frac{n}{n-d}\tau^2 \epsilon^2\right).
\end{eqnarray}
By the definition of $\epsilon$ (in terms of $\delta$),
$\beta \ge 1 - \xi - \delta + O(\tau)$,
and because this holds for every $\tau$, the result follows.

\end{proof}

\begin{proof}[Proof of Lemma \ref{lemma::mod}]
Let $f,g \in A$
be such that $||f-g||_p \le \epsilon$.
Then,
\begin{eqnarray}
\P_g{ L \le  f \le U } \hskip -0.75in &&\nonumber\\
  &=& \P_f{ L \le  f \le U } + \P_g{ L \le  f \le U } - \P_f{ L \le  f \le U } \\
  &\ge& \P_f{ L \le  f \le U } - d_{\srm TV}(P_f,P_g) \\ 
  &\ge& 1 - \alpha - M_p(||f-g||_p,A) \\
  &\ge& 1 - \alpha - M_p(\epsilon(f,p),A).
\end{eqnarray}
We also have that $\P_g{L \le g \le U} \ge 1 - \alpha$.
Hence,
\begin{eqnarray}
\P_g{L \le g \le U, L \le f \le U} \hskip -0.75in &&\nonumber\\
   &\ge& \P_g{ L \le g \le U } + \P_g{ L \le  f \le U } - 1 \\
   &\ge& 1 - \alpha + 1 - \alpha - M_p(\epsilon(f,p),A)  - 1 \\
   &\ge& 1 - 2\alpha - M_p(\epsilon(f,p),A).
\end{eqnarray}
The event
$\Event{L \le g \le U, L \le f \le U}$
implies that $W \ge \norm{g- f}_\infty$.
Hence,
\begin{eqnarray*}
\P_f{W > ||f-g||_\infty} &\ge & 1 - 2\alpha - M_p(\epsilon(f,p),A) \\
  &\ge& 1 - 2\alpha -  M_p(\epsilon(f,p),A) \\
  &\ge& 1 - 2\alpha -  M_p(\epsilon,A).
\end{eqnarray*}
It follows then that
\begin{equation}
\P_f{W > \epsilon(f,\infty)} = \inf_g \P_f{W > ||f-g||_\infty}.
\end{equation}
and thus
\begin{equation}
\inf_{f\in A_0} \P_f{W > \epsilon(f,\infty)} \ge 
1 - 2\alpha -  \sup_{f\in A_0}M_p(\epsilon(f,p),A).
\end{equation}
This proves the first claim.
But 
$\epsilon(f,\infty) \ge \epsilon(f,p)$ for any $1 \le p \le \infty$.
The final claim follows immediately.
\end{proof}

\begin{proof}[Proof of Lemma \ref{lemma::m2}]
Choose $f\in A_0$.
Choose $g\in A_1$ 
to minimize
$d_{\rm TV}(p_f,p_g)$ such to 
such that
$||f-g||_\infty = \epsilon$.
Hence,
$d_{\rm TV}(p_f,p_g) = m_\infty(\epsilon,A_0,A_1)$.
Then,
\begin{eqnarray}
\P_f{ L \le  g \le U } \hskip -0.75in &&\nonumber\\
  &=& \P_g{ L \le  g \le U } + \P_f{ L \le  g \le U } - \P_g{ L \le  g \le U } \\
  &\ge& \P_g{ L \le  g \le U } - d_{\srm TV}(P_f,P_g) \\  
  &\ge& 1 - \alpha - m_\infty(\epsilon,A_0,A_1)
\end{eqnarray}
because, by assumption. $\P_g{L \le g \le U} \ge 1 - \alpha$.
We also have that $\P_f{L \le f \le U} \ge 1 - \alpha$.
Hence,
\begin{eqnarray}
\P_f{L \le f \le U, L \le g \le U} \hskip -0.75in &&\nonumber\\
   &\ge& \P_f{ L \le f \le U } + \P_f{ L \le  g \le U } - 1 \\
   &\ge& 1 - \alpha + 1 - \alpha - m_\infty(\epsilon,A_0,A_1)\\
   &\ge& 1 - 2\alpha - m_\infty(\epsilon,A_0,A_1).
\end{eqnarray}
The event
$\Event{L \le f \le U, L \le g \le U}$
implies that $W \ge \norm{f- g}_\infty$.
Hence,
\begin{equation}
\P_f{W > ||f-g||_\infty} \ge  1 - 2\alpha - m_\infty(\epsilon,A_0,A_1).
\end{equation}
It follows then that
\begin{equation}
\sup_{f\in A_0}\P_f{W > \epsilon} \ge 1 - 2\alpha - m_\infty(\epsilon,A_0,A_1).
\end{equation}
\end{proof}

\begin{proof}[Proof of Theorem \ref{thm::target-rate}]
First, we compute $m_\infty(\epsilon,\cF,\cF)$.
Note that for all $f\in\cF$, $d_{TV}(\P_f,\P_0) = \tau(\sqrt{n}\norm{f})$.
Hence, 
$m_\infty(\epsilon,\cF,\cF) = \tau(\sqrt{n} v)$
where 
$v = \min\{ ||f||:\ f\in\cF,\ \norm{f}_\infty = \epsilon \}$.
By Lemma \ref{lemma::min2inInfinity}, 
$v=\epsilon/(\sqrt{n}\Omega_\cF)$.
It follows by Lemma \ref{lemma::m2} that
\begin{equation}
\sup_{f\in\cF}\P{W > w} \ge  1 - 2\alpha - \tau\left(\frac{w}{\sigma\Omega_\cF}\right).
\end{equation}
Let $w_*= \sigma\Omega \tau^{-1}(1-2\alpha-\gamma)$.
It follows that
if $w < w_*$ then
$\inf_{f\in\cF}\P{W \le w} <  1 - \gamma$ 
which is a contradiction.

That the proposed band has correct coverage follows easily.
Now, $(\Pi \Pi^T)_{ii} \le \Omega_{\cF}$ and
$z_{\alpha/2n} \le \sqrt{c\log n}$ for some $c$
and the claim follows.
\end{proof}

\begin{proof}[Proof of Theorem \ref{thm::main}]
We break the argument up into three parts.
Parts I and II taken together contribute the term $v_0$ from equation (\ref{eq::v0const}) to the bounds.
The logic of both parts is the same: find a value $w_*$ such that if $w < w_*$ then $\sup_{f\in\cF} \P{W > w} > \gamma$.
and, equivalently, $\inf_{f\in\cF}\P{W \le w} < 1 - \gamma$, which gives a contradiction under the assumptions of the theorem.
Part III contributes the term $v_1$ from equation (\ref{eq::v1const}) to the bounds.
It is based on using the confidence bands to construct both an estimator and a test.
Throughout the proof, we refer to the space $V \supset \cF$ defined in equation (\ref{eq::V}); this 
is the set of spoilers that are within $\epsilon_2$ of $\cF$.

\smallskip

{\sf Part I.}
First, we compute $m_\infty(w,\cF,\cF)$.
Note that for all $f\in\cF$, $d_{\srm TV}(\P_f,\P_0) = \tau(\sqrt{n}\norm{f}/\sigma)$.
Hence, 
$m_\infty(w,\cF,\cF) = \tau(\sqrt{n} v/\sigma)$
where 
$v = \min\{ ||f||:\ f\in\cF,\ \norm{f}_\infty = \epsilon \}$.
By Lemma \ref{lemma::min2inInfinity}, 
$v = w/(\sqrt{n}\Omega_\cF)$.
It follows by Lemma \ref{lemma::m2} that
\begin{equation}
\sup_{f\in\cF}\P{W > w} \ge  1 - 2\alpha - \tau\left(\frac{w}{\sigma\Omega_\cF}\right).
\end{equation}
Take $w_*= \sigma\Omega_\cF \tau^{-1}(1-2\alpha-\gamma)$.

\medskip

{\sf Part II.} 
\emph{Case (a.)}\enspace 
$\ds \epsilon_2 \le \epsilon_\infty/\sqrt{n}$.
\enspace 
First, note that 
$m_\infty(w,\cF,V) = \tau(\sqrt{n}\frac{w}{\sigma\sqrt{n}}) = \tau(w/\sigma)$
for $w \le \sqrt{n}\epsilon_2$,
because the minimum two-norm for a given infinity-norm is achieved
on the coordinate axis.
Second, let $A_0=\cF$ and $A_1 = V$ in Lemma \ref{lemma::m2}.
Then, for $w \le \sqrt{n}\epsilon_2$,
\begin{equation}
\sup_{f\in\cF}\P{W > w} \ge  1 - 2\alpha - \tau\left(\frac{w}{\sigma}\right)
\end{equation}
Let $w_* = \sigma\min(\tau^{-1}(1 - 2\alpha - \gamma),\epsilon_2\sqrt{n})$, 
then
$\sup_{f\in\cF}\P{W > w_0} \ge  \gamma.$

\emph{Case (b.)}\enspace 
$\ds \epsilon_2 > \epsilon_\infty/\sqrt{n}$. \enspace 
First, note that
$m_\infty(w,\cF,V) = \tau(\sqrt{n}\frac{w}{\sigma\sqrt{n}}) = \tau(w/\sigma)$
for $w \le \epsilon_\infty$.
Second, let $A_0=\cF$ and $A_1 = V$ in Lemma \ref{lemma::m2}.
Then, for $w \le \epsilon_\infty$,
\begin{equation}
\sup_{f\in\cF}\P{W > w} \ge  1 - 2\alpha - \tau\left(\frac{w}{\sigma}\right)
\end{equation}
Let $w_* = \sigma\min(\tau^{-1}(1 - 2\alpha - \gamma),\epsilon_\infty)$, 
then
$\sup_{f\in\cF}\P{W > w_0} \ge  \gamma.$

{\sf Part III.} 
The argument here is based on an argument in Baraud (2004).
Let $\hat f = (U + L)/2$.
Define a rejection region
\begin{equation}
\cR = \Set{ W > w} \cup \Set{ ||\hat{f} - \Pi \hat{f}||_2 > \frac{W}{2} }.
\end{equation}
Now, for any $f\in \cF$, $f^\star=f$,
$||\hat{f} - \Pi \hat{f}||_2  \le ||\hat{f} - f||_2$ and
\begin{eqnarray}
\P_f(\cR)
  &\le& \P_f{ W > w} + \P_f{ ||\hat{f} - \Pi \hat{f}||_2 > W/2}\\
  &\le& \gamma + \P_f{ ||\hat{f} - \Pi \hat{f}||_2 > W/2}\\
  &\le& \gamma + \P_f{ ||f - \hat{f}||_2 > W/2}\\
  &=&   \gamma + \P_f{ ||f^\star - \hat{f}||_2 > W/2}\\
  &\le& \gamma + \P_f{ ||f^\star - \hat{f}||_\infty > W/2}\\
  &\le& \gamma + \alpha
\end{eqnarray}
which bounds the type I error of $\cR$.

Now let $f$ be such that
$\norm{f - \Pi f}  > \max\{w,\epsilon_2\}$.
Because $\norm{f - \Pi \hat f} > \norm{f - \Pi f}$,
$\norm{f - \Pi f}  > \epsilon_2$ implies that $f^\star = f$.
And thus,
\begin{equation}
||\hat{f} - \Pi\hat{f}||_2 
  \ge ||f - \Pi\hat{f}||_2 - ||f-\hat{f}||_2 
  \ge w - ||f-\hat{f}||_2.
\end{equation}
Hence,
\begin{eqnarray}
\P_f(\cR^c) 
   &=&   \P_f{  ||\hat{f} - \Pi \hat{f}||_2 \le W/2, W/2\le w/2}\\
   &\le& \P_f{  ||\hat{f} - \Pi \hat{f}||_2 \le w/2, W\le w}\\
   &\le& \P_f{  ||f-\hat{f}||_2 \ge w/2, w\ge W}\\
   &\le& \P_f{  ||f-\hat{f}||_2 \ge W/2}\\
   &=&   \P_f{  ||f^\star-\hat{f}||_2 \ge W/2}\\
   &\le& \P_f{  ||f^\star-\hat{f}||_\infty \ge W/2}\\
   &\le& \alpha.
\end{eqnarray}
Thus, $\cR$ defines a test 
for $H_0:f\in\cF$
with level $\alpha+\gamma$
whose power 
more than a distance
$\max\{w,\epsilon_2\}$ from $\cF$ is at least $1-\alpha$.
Using Lemma \ref{lemma::baraud-doover} with $\xi = \alpha + \gamma$
and $\delta = 1 - \gamma - 2\alpha$, this implies that
\begin{equation}
\max\{ w,\epsilon_2\} \ge 2 \kappa(\alpha,\gamma) (n - d)^{1/4} n^{-1/2}.
\end{equation}
The result follows.
\end{proof}

\begin{proof}[Proof of Theorem \ref{thm::optimal}]
The volume is minimized by making $\epsilon_\infty$ as large as possible
and
$\epsilon_2$ as small as possible.
To achieve the lower bound on the width
requires
$\epsilon_\infty \le w_\cF$ and
$\epsilon_2 \ge 2\kappa(\alpha,\gamma)(n-d)^{1/4}n^{-1/2}.$
\end{proof}

\begin{proof}[Proof of Theorem \ref{thm::achieve}]
Let $A = \Event{T \le \chi^2_{\gamma,n-d}}$.
Then,
$$
\P_f{f^\star\notin B} = \P_f{f^\star\notin B,A} + \P_f{f^\star\notin B,A^c}.
$$
We claim that
$\P_f{f^\star\notin B,A}\le \alpha/2$
and 
$\P_f{f^\star\notin B,A^c} \le \alpha/2$.
There are four cases.

\smallskip
{\em Case I.} $f\in\cF$. Then $f=f^\star$ and
$\P_f{f\notin B,A^c} \le \P_f{A^c} \le \alpha/2$.
$\P_f{f\notin B,A}\le \P_f{f\notin B} =
\P_{\Pi f}{\Pi f\notin B} \le \P_{\Pi f}{||\hat f - \Pi f||_\infty > w_{\cF}} \le \alpha/2$.

\bigskip

{\em Case II.} $f\in V-\cF$ where
$V= \{f:\ \norm{f-\Pi f}\le \epsilon_2, \, \norm{f-\Pi f}_\infty \le \epsilon_\epsilon\}$.
Again, $f=f^\star$.
First,
$\P_f{f\notin B,A^c} \le  \P_f{||Y-f||_\infty > z_{\alpha/2n}} \le  \alpha/2$.
Next, we bound
$\P_f{f\notin B,A}$.
Note that
$\hat{f} = \Pi Y \sim N(g, \sigma^2 \Pi \Pi^T)$, where $g=\Pi f$.
Then
$\hat{f}_i \sim N(g_i, \Omega_i^2)$.
Let $B_0 = (L+\epsilon_\infty, U-\epsilon_\infty)$.
Then, $\Pi f\in B_0$ implies $f\in B$ and
$\P_f{\notin B, A}\le \P_f{\Pi f \notin B_0}\le \alpha/2$.

\bigskip

{\em Case III.} $f\notin V$,
$||f-\Pi f|| \le \epsilon_2$ and
$||f-\Pi f||_\infty > \epsilon_\infty$.
In this case, $f^\star = \Pi f$.
Then
$\P_f{f^\star,f\in B^c , A^c}\le
\P_f{f\in B^c , A^c}\le \alpha/2$.
Also,
$\P_f{f^\star,f\in B^c,A}\le \P_f{f^\star\notin B} =
\P_{\Pi f}{\Pi f\notin B} \le \P_{\Pi f}{||\hat f - \Pi f||_{\infty} > w_{\cF}} \le \alpha/2$.

\bigskip

{\em Case IV.} $f\notin V$ and $||f-\Pi f|| > \epsilon_2$.
In this case, $f^\star = f$.
But
$$
\P_f{f\notin B,A}\le \P_f{A}\le 
F_{f-\Pi f,n-d}(\chi^2_{\gamma,n-d}) \le 
F_{\epsilon_2,n-d}(\chi^2_{\gamma,n-d}) \le 
\alpha/2
$$
and
$$
\P_f{f\notin B,A^c} \le \P_f{f\notin B,A^c} \le \alpha/2.
$$

\bigskip
Thus, $\P_f{f^\star \not\in B} \le \alpha$.
Equation (\ref{eq::it-works}) follows since
$\P_f{T \le \chi^2_{\gamma,n-d}} \ge 1-\gamma$ for all $f\in \cF$.
\end{proof}

\begin{proof}[Proof of Lemma \ref{lemma::Modulus-of-Continuity}]
First note that if $B$ is a ball in $\mathbb{R}^n$ in any norm,
then $B - B = 2 B$.
Second, we have that
\begin{eqnarray}
\omega(u) 
   &=& \sup\{ |T g| :\; \norm{g}_2 \le u,\  g\in V-V\} \\
   &=& \sup\{ |T g| :\; \norm{g}_2 \le u,\  g\in V(2\epsilon_2,2\epsilon_\infty)\}.
\end{eqnarray}
To see the latter equality, note that if $g,h\in V$, then we can
write $g - h = f + \delta_1 - \delta_2$ where $f \in \cF$ and $\delta_i$ are
in $B^\perp_k(0,\epsilon_k)$ for $k = 2,\infty$. Thus, $\delta_1 - \delta_2$
is in $2B^\perp_2(0,\epsilon_2) \cap 2B^\perp_\infty(0,\epsilon_\infty)$.

Set $B^*(f) = B^\perp_2(f,2\epsilon_2) \cap B^\perp_\infty(f,2\epsilon_\infty)$.
We have that
\begin{eqnarray}
\omega(\eta,\cF) &=& \sup\{ f_1:\; \norm{f}_2 \le \eta, f\in \cF\} \\
\omega(\eta, B^*(0)) &=& \sup\{ f_1:\; \norm{f}_2 \le \eta, f\in B^*(0)\}.
\end{eqnarray}
For any $g\in V(2\epsilon_2,2\epsilon_\infty)$, we can write $g = g_1 + g_2$
where $g_1 \in \cF$ and $g_2 \in B^*(0)$ and the two functions are orthogonal.
Then,
\begin{eqnarray}
w(u,V) 
   &=& \sup \Biggl\{T(g):\ g\in V(2\epsilon_2,2\epsilon_\infty),\  \norm{g}_2 \le u\Biggr\}\\
   &=& \sup_{0\le c \le u} \Biggl\{T(g_1 + g_2):\ \norm{g_1}_2 \le \sqrt{u^2 - c^2}, \,\norm{g_2}_2 \le c^2,\Biggr. \nonumber\\
   && \phantom{\sup_{0\le c \le u} \Biggl\{T(g_1 + g_2):\ \Biggr.}
          \   g_1\in \cF, g_2 \in B^*(0)\Biggr\}\\
   &\le& \sup_{0 \le c \le u} \left[ \sup_{g_1\in\cF\atop\norm{g_1}_2 \le \sqrt{u^2 - c^2}} T(g_1) + 
                \sup_{g_2 \in B^*(0)\atop \norm{g_2}_2\le c} T(g_2) \right] \\
   &=&  \sup_{0 \le c \le u} \left[\omega(\sqrt{u^2 - c^2},\cF) + \omega(c,B^*(0))\right].
\end{eqnarray}
Moreover, equality can be attained for each $c$ by choosing $g_1$ and $g_2$ to be the maximizers
(or suitably close approximants thereof) of each term in the last equation.
Consequently,
\begin{equation}
\omega(u) = \sup_{0 \le c \le u} \omega(\sqrt{u^2 - c^2},\cF) + \omega(c,B^*(0)).
\end{equation}

To derive $\omega(\eta,B^*(0))$, note that 
$f = ((\eta\wedge \epsilon_2)\sqrt{n} \wedge\epsilon_\infty,0,0,\ldots,0)$
maximizes $f_1$ subject to the norm constraint.
Hence, $\omega(\eta,B^*(0)) = \min( (\eta\wedge \epsilon_2)\sqrt{n}, \epsilon_\infty )$.
For $\omega(\eta,\cF)$, let $e = (1,0,\ldots,0) \in\mathbb{R}^n$.
Recall that 
$\Omega_{\cF} = \frac{\langle e, \Pi_\cF e\rangle}{\norm{e}\,\norm{\Pi_\cF e}} = \frac{\norm{\Pi_\cF e}}{\norm{e}}$,
which is between 0 and 1.
Maximizing $e^T f$ for $f\in\cF$ and $\norm{f}_2\le \eta$
is equivalent to maximizing $n \langle e, f\rangle = n \langle \Pi_\cF e, f\rangle$.
The maximum subject to the constraint occurs at $f^\star = \eta \Pi e/\norm{\Pi e}$
Hence, $\omega(\eta,\cF) = \eta\sqrt{n} \Omega_{\cF}$.
Note that $\eta$ is in terms of the normalized two norm;
in the ``natural'' (root sum of squares) norm, the modulus would be $\omega_\natural(u,\cF) = u \Omega_{\cF}$.

It follows that
\begin{eqnarray}
\omega(u,V) \hskip -0.25in &&\nonumber\\
  &=& \sup_{0 \le c \le u} [\omega(\sqrt{u^2 - c^2},\cF) + \omega(c,B^*(0))] \\
  &=& \sup_{0 \le c \le u} \left[\sqrt{n}\Omega_{\cF}\sqrt{u^2 - c^2} + \min( (c\wedge \epsilon_2)\sqrt{n}, \epsilon_\infty )\right] \\
  &=& \sqrt{n}\sup_{0 \le c \le u} \left[\Omega_{\cF}\sqrt{u^2 - c^2} + \min( c, \epsilon_2\wedge (\epsilon_\infty/\sqrt{n}))\right] \\
  &=& \sqrt{n} \left( u \Omega \sqrt{\frac{\Omega^2}{1 + \Omega^2}} + \min( \frac{u}{\sqrt{1+\Omega^2}}, \epsilon_2\wedge (\epsilon_\infty/\sqrt{n}))\right) \\
  &=& \left( u\sqrt{n} \Omega \sqrt{\frac{\Omega^2}{1 + \Omega^2}} + \min\left( \frac{u\sqrt{n}}{\sqrt{1+\Omega^2}},\, \epsilon_2\sqrt{n},\, \epsilon_\infty\right)\right)
\end{eqnarray}
because the supremum over $c$ is maximized at $c = u/(1 + \Omega^2)$.
In the natural two norm, we have
\begin{equation}
\omega_\natural(u,V) = 
\left( u \Omega \sqrt{\frac{\Omega^2}{1 + \Omega^2}} + 
\min\left( \frac{u}{\Omega}\,\sqrt{\frac{\Omega^2}{1+\Omega^2}},\, 
\epsilon_{2,\natural},\, \epsilon_\infty\right)\right).
\end{equation}
\end{proof}

Next, we prove the lower bound result generalized to a nested sequence of
subspaces. To do so, we need to prove several auxilliary lemmas.
Define for each $1 \le j \le m$,
\begin{equation}\label{eq::Uj}
U_j = \Set{ f \in \R^n:\ F^*(f) = \{\Pi_j f, f\} \ {\rm or}\ F^*(f) = \{f\} }.
\end{equation}
Referring to the definition of $V$ in equation (\ref{eq::V}), define here
$V_j = V(\cF_j,\epsilon_{2,j},\epsilon_{\infty,j})$.

\begin{lemma}\label{lemma::modulusUj}
Let $w > 0$. Then,
\begin{eqnarray}
m_\infty( w, \cF_j \cap U_j, \cF_j \cap U_j ) &=& m_\infty( w, \cF_j, \cF_j ) \\
m_\infty( w, \cF_j \cap U_j,   V_j \cap U_j ) &=& m_\infty( w, \cF_j, V_j )
\end{eqnarray}

\end{lemma}

\begin{proof}%[ of Lemma \ref{lemma::modulusUj}]

First, let $f, g \in \cF_j$ be the minimal pair for $m_\infty( w, \cF_j, \cF_j )$.
Let $\psi$ be a unit-2-norm vector in $\cF_{j} \cap \cF_{j-1}^\perp$.
Let $\lambda > \epsilon_{2,1}$ and 
define
\begin{eqnarray}
\tilde f &=& \lambda \psi + f \\
\tilde g &=& \lambda \psi + g.
\end{eqnarray}
Then, $\tilde f, \tilde g \in \cF_j \cap U_j$ because
if either $f$ or $g$ were in $\cF_j \cap U_j^c$ then adding $\lambda\psi$ makes the
distance from the projection on one of the lower spaces larger than the corresponding $\epsilon_2$.
Also $d_{\srm TV}(P_{\tilde f}, P_{\tilde g}) = d_{\srm TV}(P_f,P_g)$
and $\norm{\tilde f - \tilde g}_\infty = \norm{f - g}_\infty$.
Hence, $m_\infty( w, \cF_j \cap U_j, \cF_j \cap U_j ) \le m_\infty( w, \cF_j, \cF_j )$.
But $\cF_j \cap U_j \subset \cF_j$, so 
$m_\infty( w, \cF_j \cap U_j, \cF_j \cap U_j ) = m_\infty( w, \cF_j, \cF_j )$ as was to be proved.

Second,let $f \in \cF_j$ and $g\in V_j$ be the minimal pair for $m_\infty( w, \cF_j, V_j )$.
Now apply the same argument.

\end{proof}

\begin{lemma}\label{lemma::baraudUj}

Let $0 < \delta < 1 - \xi$ and 
\begin{equation}
\epsilon = \frac{(n-d_j)^{1/4}}{\sqrt{n}}\,\left(2 \log(1 + 4 \delta^2)\right)^{1/4}.  
\end{equation}
Define $A_j = U_j \cap \Set{f:\ \norm{f - \Pi_j f} > \epsilon}$.
Then,
\begin{equation}
\beta \equiv
   \inf_{\phi_\alpha\in\Phi_\xi}\sup_{f\in A_j}\P_f{\phi_\xi=0} \ge 1 - \xi - \delta
\end{equation}
where
\begin{equation}
\Phi_\xi = \Biggl\{ \phi_\xi:\ \sup_{f\in\cF_j}\P_f{\phi_\xi = 0} \le \xi \Biggr\}
\end{equation}
is the set of level $\xi$ tests.

\end{lemma}

\begin{proof}%[ of Lemma \ref{lemma::baraudUj}]

Let $f_E$ be defined as in equation (\ref{eq::fE})
in the proof of Lemma \ref{lemma::baraud-doover}.
Let $\psi$ be a unit vector in $\cF_{j+1} \cap \cF_j^\perp$
and let $\lambda > \epsilon_{2,1}$.
Then, define $\tilde f_E = \lambda \psi + f_E$.
Now apply the proof of Lemma \ref{lemma::baraud-doover} using $f_0 = \lambda\psi$ instead of $0$.
The total variation distances among corners of the hypercube do not change and the result follows.

\end{proof}

\begin{lemma}\label{lemma::main::special}
Fix $0< \alpha < 1$ and
$0 < \gamma < 1-2\alpha$.
Suppose that for bands $B = (L,U)$
\begin{equation} \label{eq::validband::special}
\inf_{f\in U_j} \P_f{ F^*(f) \cap B \ne \emptyset }\ge 1-\alpha. 
\end{equation}
Then 
\begin{equation}\label{eq::Wgamma::special}
\inf_{f\in\cF_j} \P_f{W \le w} \ge 1 - \gamma.
\end{equation}
implies
\begin{equation}\label{eq::the-cool-lower-bound::special}
w \ge \underline{w}(\cF_j,\epsilon_{2,j},\epsilon_{\infty,j},n,d_j,\alpha,\gamma,\sigma),
\end{equation}
where $\underline{w}$ is given in Theorem \ref{thm::main}.
\end{lemma}

\begin{proof}%[ of Lemma \ref{lemma::main::special}]

To prove this lemma, we will adapt the proof of Theorem \ref{thm::main} as follows.
By Lemma \ref{lemma::modulusUj},
the argument for Parts I and II is the same with $\cF$ replaced with $\cF_j \cap U_j$ 
and $V$ replaced with $V_j \cap U_j$.
By replacing the reference to Lemma \ref{lemma::baraud-doover} with Lemma \ref{lemma::baraudUj}, 
the argument for Part III also follows exactly.
The result follows.

\end{proof}

\begin{proof}[Proof of Theorem \ref{thm::main::nested}]

The result follows directly from Lemma \ref{lemma::main::special} because
$\inf_{f\in\R^n} \P{ F^*(f) \cap B \ne \emptyset} \ge 1 - \alpha$ implies 
$\inf_{f\in U_j} \P{ F^*(f) \cap B \ne \emptyset} \ge 1 - \alpha$.

\end{proof}

\begin{proof}[Proof of Theorem \ref{thm::achieve-multi}]
Note that
$\P_f{ F^\star \cap B = \emptyset} = \sum_j \P_f{ F^\star \cap B = \emptyset,\hat{J}=j}$.
We show that
$\P_f{ F^\star \cap B = \emptyset,\hat{J}=j} \le \alpha_j$ for each $j$.
There are three cases. Throughout the proof, we take $\sigma = 1$.

\emph{Case I.} $||f-\Pi_j f|| > \epsilon_{2,j}$.
Then,
\begin{eqnarray*}
\P_f{ F^\star \cap B = \emptyset,\hat{J}=j}   &\le&
\P_f{\hat{J}=j} \le F_{f-\Pi_j f,n-d_j}(\chi^2_{\gamma,n-d_j})\\
  &\le&F_{\epsilon_{2,j},n-d_j}(\chi^2_{\gamma,n-d_j})\\
  &\le& \alpha_j
\end{eqnarray*}
due to (\ref{eq::gam-alph}).

\emph{Case II.} $||f-\Pi_j f|| \le \epsilon_{2,j}$ and
$||f-\Pi_j f||_\infty \le \epsilon_{\infty,j}$.
So,
\begin{eqnarray*}
\P_f{ F^\star \cap B = \emptyset,\hat{J}=j}   &\le& \P_f{ f \notin B,\hat{J}=j}\\
  &\le& \P_f{ ||f-\hat{f}||_\infty > w_{\cF_j} + \epsilon_{\infty,j}}\\
  &\le& \P_f{ ||f-\Pi_j f||_\infty +||\Pi_j f-\Pi_j Y||_\infty > w_{\cF_j} + \epsilon_{\infty,j}}\\
  &\le& \P_f{ ||\Pi_j f-\Pi_j Y||_\infty > w_{\cF_j}}\\
  &=& \P_{\Pi_j f}{ ||\Pi_j f-\Pi_j Y||_\infty > w_{\cF_j}}\\
  &\le& \alpha_j.
\end{eqnarray*}

\emph{Case III.}
$||f-\Pi_j f|| \le \epsilon_{2,j}$ and
$||f-\Pi_j f||_\infty > \epsilon_{\infty,j}$.
Now,
\begin{eqnarray*}
\P_f{ F^\star \cap B = \emptyset,\hat{J}=j} 
  &\le& \P_f{ \Pi_j f \notin B,\hat{J}=j}\\
    &=& \P_f{ \norm{\Pi_j Y - \Pi_j f}_\infty > c_j,\hat{J}=j}\\
  &\le& \P_f{ \norm{\Pi_j Y - \Pi_j f}_\infty > c_j}\\
    &=& \P_{\Pi_j f}{\norm{\Pi_j Y - \Pi_j f}_\infty > c_j}\\
  &\le& \alpha_j. 
\end{eqnarray*}

To prove (\ref{eq::the-gam-result}),
suppose that $f\in\cF_j$.
Then,
$\P_f{\hat{J} > j} \le \gamma$.
But, as long as $\hat{J}\le j$,
$W = w_{\cal_{\hat J}}(\alpha_{\hat J}) + \epsilon_{\infty,\hat{J}} \le
w_j(\alpha_{j}) + \epsilon_{\infty,j}$.
The last statement follows since,
when $\epsilon_{2,j} \ge Q(n-d_j,\alpha/2,\gamma)(n-d_j)^{1/4} n^{-1/2}$
\end{proof}

\section{Discussion}\label{sec::discussion}

We have shown that adaptive confidence bands for
$f$ are possible if coverage is replaced by 
surrogate coverage.
Of course, there are many other ways one could
define a surrogate.
Here, we briefly outline a few possibilities.

Wavelet expansions of the form
$$
f(x) = \sum_j \alpha_j \phi_j(x) + \sum_j\sum_k \beta_{jk}\psi_{jk}
$$
lend themselves quite naturally to the surrogate approach.
For example, one can define
$$
f^\star(x) = \sum_j \alpha_j \phi_j(x) + \sum_j\sum_k s(\beta_{jk})\psi_{jk}
$$
where $s(x) = {\rm sign}(x)(|x|-\lambda)_+$ is the usual soft-thresholding function.

For kernel smoothers and local polynomial smoothers $\hat{f}_h$ that
depends on a bandwidth $h$, a possible surrogate is
$f^\star = \E(\hat{f}_{h^\star})$
where $h^\star$ is the largest bandwidth $h$ for which
$\hat{f}_h$ passes a goodness of fit test with high probability.
In the spirit of Davies and Kovac (2001),
one could take the test to be a test for randomness
applied to the residuals.

Motivated by ideas in Donoho (1988)
we can define another surrogate as follows.
Let us switch to the problem of
density estimation.
Let $X_1, \ldots, X_n \sim F$
for some distribution $F$.
The goal is define an appropriate surrogate band for the density $f$.
Define the smoothness functional
$S(F) = \int (f''(x))^2 dx$.
To make sure that $S(F)$ is well defined for all $F$
we borrow an idea from Donoho (1988).
Let $\Phi_h$ denote a Gaussian with standard deviation $h$ and define
$S(F) = \lim_{h\to 0}S(F\oplus \Phi_h)$
where $\oplus$ denote convolution.
Donoho shows that $S$ is then a well-defined, convex, lower semicontinuous functional.

Let $\hat{F}_n$ be the empirical distribution function and let
$B=B(\hat{F},\epsilon_n) = \{ F:\ ||F - \hat{F}_n|| \le \epsilon_n\}$
where
$||\cdot ||$ is the Kolmogorov-Smnirnov distance and
$\epsilon_n$ is the $1-\beta$ quantile of 
$||U-U_n||$ where $U$ is the uniform distribution
and $U_n$ is the empirical from a sample from $U$.
Thus, $B$ is a nonparametric, $1-\beta$ confidence ball for $F$.
The simplest $F\in B$ is the distribution
that minimize $S(F)$ subject to
$F\in B$.
We define the surrogate $F^\star$ to be the distribution that
minimizes $S(F)$ subject to $F$ belonging to $B_F$, 
where $B_F$ is a population
version of $B$.
We might then think of $F^\star$ as the simplest distribution
that is not empirically dinstinguishable from $F$.
A natural definition of $B_F$ might be
$B_F = \{ G:\ ||F-G|| \le \epsilon_n\}$.
But this definition only makes sense for fixed radius 
confidence sets.
Another definition is
$B_F = \{ G :\ \P_F{G\in B}\ge 1/2\}$.

To summarize, we define 
\begin{equation}
F^\star = {\rm argmin}_{F \in B_F} S(F)
\end{equation}
where
\begin{equation}
B_F = \Biggl\{ G :\ \P_F{G\in B(\hat{F}_n,\epsilon_n)}\ge 1/2\Biggr\}
\end{equation}
and
$B(\hat{F}_n,\epsilon_n) = \{ G:\ ||\hat{F}_n-G|| \le \epsilon_n\}$.
Let
\begin{equation}
\Gamma = \cup \{G^\star:\ G\in B(\hat{F}_n,\epsilon_n)\}.
\end{equation}
Then
\begin{equation}
\ell(x) = \inf_{F\in\Gamma}F'(x),\ \ \ 
u(x) = \sup_{F\in\Gamma}F'(x)
\end{equation}
defines a valid confidence band for the density of $F^\star$.

Let us also mention
average coverage (Wahba 1983; Cummins, Filloon, Nychka 2001).
Bands $(L,U)$ have average coverage if
$\P_f{L(\xi) \le f(\xi) \le U(\xi)} \ge 1-\alpha$
where $\xi\sim {\rm Uniform}(0,1)$.
A way to combine average with the surrogate idea is to enforce
something stronger than average coverage such as
$$
\P_f{L(\xi) \le f(\xi) \le U(\xi)\ {\rm and}\  \hat{f}\preceq f} \ge 1-\alpha
$$
where $\hat{f} = (L+U)/2$ and
$\hat{f}\preceq f$ means that
$\hat{f}$ is simpler than $f$ according to a
partial order $\preceq$, for example,
$f\preceq g$ if
$\int (f'')^2 \leq \int (g'')^2$.

\section*{References}\hfil\break

\indent Baraud, Y. (2002).
Non Asymptotic minimax rates of testing in signal detection,
\emph{Bernoulli}, {\bf 8}, 577.

Baraud, Y. (2004).
Confidence balls in {G}aussian regression,
{\em The Annals of Statistics}, 32, 528--551.

Beran, Rudolf and D\"{u}mbgen, Lutz. (1998).
Modulation of estimators and confidence sets.
{\em The Annals of Statistics}, 26, 1826--1856.

Bickel, P.J. and Ritov, Y. (2000).
Non-and semi parametric statistics: compared and contrasted.
{\em J. Statist. Plann. Inference}, 91, 

Birg\'e, L. (2001).
An alternative point of view on Lepski's method.
In {\em State of the Art in Probability and Statistics.}
(M. de Gunst, C. Klaassen and A. van der Vaart, eds.)
113--133, IMS, Beachwood, OH.

Cai, T. and Low, M. (2005).
Adaptive Confidence Balls.
{\em The Annals of Statistics}, 34, 202--228.

Cai, T. and Low, Mark, G. (2004).
An adaptation theory for nonparametric confidence intervals.
{\em Ann. Statist.},  32, 1805--1840.

Chaudhuri, Probal and Marron, J. S. (2000).
Scale space view of curve estimation.
{\em The Annals of Statistics}, 28, 408--428.

Claeskens, G. and Van Keilegom, I. (2003).
Bootstrap confidence bands for regression curves
and their derivatives.
{\em The Annals of Statistics}, 31, 1852--1884.

Cummins D., Filloon T., Nychka D. (2001).
Confidence Intervals for Nonparametric Curve Estimates: 
Toward More Uniform Pointwise Coverage
{\em Journal of the American Statistical Association}, 96, 233--246.

Donoho, D. (1988).
One-Sided Inference about Functionals of a Density.
{\em Annals of Statistics}, 16, 1390--1420.

Donoho, D. (1995).
De-noising by soft-thresholding.
{\em IEEE Transactions on Information Theory}, 41, 613--627.

Donoho, D. and Liu, R. (1991).
Geometrizing Rates of Convergence, II. 
{\em The Annals of Statistics}, 19, 633--667.

Donoho, D., Johnstone, I.M., Kerkyacharian G., and Picard, D. (1995).
Wavelet Shrinkage: Asymptopia,
\emph{J. Roy. Statist. Soc. B}, {\bf 57}, 301--369.

Eubank, R.L. and Speckman, P.L. (1993).
Confidence Bands in Nonparametric Regression. 
{\em Journal of the American Statistical Association}, 
88, 1287--1301.

Genovese, C. and Wasserman, L. (2005).
Nonparametric confidence sets for wavelet regression.
{\em Annals of Statistics}, 33, 698--729.

Hall, P. and Titterington, M. (1988).
On confidence bands in nonparametric density estimation and regression.
{\em Journal of Multivariate Analysis}, 27, 228--254.

H\"{a}rdle, Wolfgang and Bowman, Adrian W. (1988).
Bootstrapping in nonparametric regression: 
Local adaptive smoothing and confidence bands.
{\em Journal of the American Statistical Association}, 83, 102--110.

H\"{a}rdle, W. and Marron, J. S. (1991).
Bootstrap simultaneous error bars for nonparametric regression.
{\em The Annals of Statistics}, 19, 778--796.

Ingster, Y. (1993).
Asymptotically minimax hypothesis testing for nonparametric alternatives, I and II.
\emph{Math. Methods Statist}, {\bf 2}, 85--114.

Ingster, Y. and Suslina, I. (2003).
{\em Nonparametric Goodness of Fit Testing Under Gaussian Models.}
Springer. New York.

Juditsky, A. and Lambert-Lacroix, S. (2003).
Nonparametric confidence set estimation.
{\em Mathematical Methods of Statistics}, 19, 410-428.

Leeb, H. and P\"otscher, B.M. (2005).
Model Selection and Inference: Facts and Fiction.
{\em Econometric Theory}, 21, 21--59.

Li, Ker-Chau. (1989).
Honest confidence regions for nonparametric regression.
{\em The Annals of Statistics}, 17, 1001--1008.

Low, Mark G. (1997).
On nonparametric confidence intervals.
{\em The Annals of Statistics}, 25, 2547--2554.

Neumann, Michael H. and Polzehl, J\"{o}rg. (1998).
Simultaneous bootstrap confidence bands in nonparametric regression.
{\em Journal of Nonparametric Statistics}, 9, 307--333.

Robins, J. and van der Vaart, Aad. (2006).
Adaptive Nonparametric Confidence Sets.
{\em The Annals of Statistics}, 34, 229--253.

Ruppert, D. and Wand, M.P. and Carroll, R.J. (2003).
{\em Semiparametric Regression},
Cambridge University Press. Cambridge.

Sun, J. and Loader, C. R. (1994).
Simultaneous confidence bands for linear regression and smoothing.
{\em The Annals of Statistics}, 22, 1328--1345.

Terrell, G.R. and Scott, D.W. (1985).
Oversmoothed Nonparametric Density Estimates.
{\em Journal of the American Statistical Association}, 
80, 209--214.

Terrell, G.R. (1990).
The Maximal Smoothing Principle in Density Estimation.
{\em Journal of the American Statistical Association}, 
85, 470--477.

Wahba, G. (1983).
Bayesian ``confidence intervals'' for the cross-validated smoothing spline.
{\em Journal of the Royal Statistical Society, Series B, Methodological}, 45, 133--150.

Xia, Y. (1998).
Bias-Corrected Confidence Bands in Nonparametric Regression. 
{\em Journal of the Royal Statistical Society. Series B}, 60, 797--811.

\end{document}